\documentclass[12pt]{amsart}
\usepackage{graphicx}
\usepackage{subfig}
\usepackage{tabularx}
\usepackage{stmaryrd}

\usepackage{mathrsfs}     % selects Times Roman as basic font
\usepackage{helvet}         % selects Helvetica as sans-serif font
\usepackage{courier}        % selects Courier as typewriter font
\usepackage{type1cm}      % activate if the above 3 fonts are
\usepackage[dvipsnames]{xcolor}
\usepackage{todonotes}
\usepackage{mathbbol}
\usepackage{multirow}
\usepackage{braket}
\usepackage{siunitx}
\usepackage{booktabs}
\usepackage{mathabx}

\usepackage{{soul}}

\usepackage{tikz}
\usepackage{pgfplots}
\usepackage{pgfplotstable}
\pgfplotsset{compat=newest}

\usepackage{geometry,calc,color}
\usepackage{amsmath,amssymb}

\usepackage{enumerate}

\renewcommand{\phi}{\varphi}
\renewcommand{\hat}[1]{\widehat{#1}}

\newtheorem{theorem}{Theorem}[section]
\newtheorem{corollary}[theorem]{Corollary}
\newtheorem{lemma}[theorem]{Lemma}
\newtheorem{proposition}[theorem]{Proposition}
\theoremstyle{definition}

\newtheorem{problem}{Problem}[section]
\newtheorem{assumption}{Assumption}[section]
\newtheorem{comment}{Comment}

\theoremstyle{remark}
\newtheorem{remark}[theorem]{Remark}
\newtheorem{property}{Property}[section]

\numberwithin{equation}{section}

\usepackage{hyperref}

\newcommand{\blu}[1]{{\color{black} #1}} % blu -> nero
\newcommand{\rosso}[1]{{\color{red} #1}} %rosso
\newcommand{\tb}{\interleave}

\newcommand{\D}{D}
\newcommand{\Tess}{\mathcal{T}_h}
\newcommand{\Edges}{\mathcal{E}_h}
\newcommand{\Squel}{\Sigma}
\renewcommand{\Squel}{\Edges}

\newcommand{\roundPrecision}{2}

\newcommand{\normal}{{\mathbf{n}}}
\newcommand{\normalK}{\normal_K}

\newcommand{\scalminusone}[2]{(#1,#2)_{-1,K}}

\newcommand{\w}{u}
\newcommand{\la}{\lambda}
\newcommand{\uh}{u_h}
\newcommand{\lah}{\lambda_h}
\newcommand{\hlah}{\widehat \lambda_h}

\date{}

\begin{document}
\title[\MakeLowercase{$p$}-robust negative norm stabilization
for DG methods]{\blu{ A \MakeLowercase{\large $p$}-robust polygonal discontinuous Galerkin method with minus one stabilization}}

\author[S. Bertoluzza]{Silvia Bertoluzza$^\sharp$} 
\author[I. Perugia]{Ilaria Perugia$^\mathsection$}
\thanks{S. Bertoluzza and I. Perugia would like to acknowledge the
  kind hospitality of the Erwin Schr\"odinger International Institute
  for Mathematics and Physics (ESI), where part of this research was
  developed under the frame of the Thematic Programme {\it Numerical
    Analysis of Complex PDE Models in the Sciences}.\\
  S. Bertoluzza and D. Prada acknowledge the support of  European Research Council (ERC), under the EU's Horizon 2020 research and innovation programmes (Project CHANGE, grant agreement No 694515).
I.~Perugia has been funded by the Austrian Science Fund (FWF) through the projects P~29197-N32 and F~65.
}
\author[D. Prada]{Daniele Prada$^\sharp$}

\date{\today}

\maketitle 

\vspace*{-0.6cm}
\begin{center}
{\small 
$^\sharp$ {IMATI ``Enrico Magenes'', CNR\\
Via Ferrata 1, 27100 Pavia, Italy\\
silvia.bertoluzza@imati.cnr.it, daniele.prada@imati.cnr.it\\ \medskip
$^\mathsection$ Faculty of Mathematics, University of Vienna\\
				Oskar-Morgenstern-Platz 1, 1090 Vienna, Austria\\
ilaria.perugia@univie.ac.at
}}
\end{center}
\vspace*{0.2cm}

% ----------------------------------------------------------------
\begin{abstract}
We introduce a new stabilization for discontinuous Galerkin methods for
  the Poisson problem on polygonal meshes, which induces optimal
  convergence rates in the polynomial approximation
  degree $p$. In the setting of [S. Bertoluzza and D. Prada, A
  polygonal discontinuous Galerkin method with minus one
  stabilization, \emph{ESAIM Math. Mod. Numer. Anal.} (DOI:
  10.1051/m2an/2020059)], the stabilization is obtained by penalizing,
  in each mesh element~$K$, a residual in the~norm of the dual of $H^{1}(K)$. This
  negative norm
  is algebraically realized via the introduction of new auxiliary
  spaces. We carry out a $p$-explicit stability and error analysis, proving
  $p$-robustness of the overall method. The theoretical findings are
  demonstrated in a series of numerical experiments.
\end{abstract}
	
\medskip\noindent
\textbf{AMS subject classification}: 
	

% ----------------------------------------------------------------
\section{Introduction}\label{sec:1}
Polytopal methods for the solution of partial differential equations
have, in recent year, gained an increased popularity thanks to the
flexibility inherently offered by the use of polytopal meshes. Indeed
polygonal meshes allow to take into account the geometrical feature of
the physical domain without resulting in an excess of degrees of
freedom, they can be used, by agglomeration, as a transition step when
dealing with triangular/tetrahedral and quadrilateral/hexahedral
meshes, and they allow for simple refining and coarsening strategies
in the framework of adaptive methods. Among the different approaches, besides Discontinuous Galerkin and its variants, such as the Hybridizable Discontinuous Galerkin method  or the Discontinuous Petrov Galerkin  method, we recall the Virtual Element method, the Hybrid High Order method  and the Mimetic Finite Differences method (see \cite{CangianietalBook}, \cite{CockburnFuSayas},  \cite{DPGII}, \cite{basicVEM},\cite{DiPietroErnLemaire} and  \cite{BrezziLipnikovSimoncini}). A key feature in all such methods is the need to resort
to some form of stabilization, which can lead to a loss of optimality
{of their convergence rates with respect to}
%both in
the mesh size~$h$ (when the mesh contains very small edges, as
compared to the element diameters) or %in
{to the polynomial degree~$p$, or both.

In~\cite{SilviaDaniele}, the focus was 
on the first issue, namely the loss of optimality with respect to the
mesh size $h$.} In that %such a
paper, a stabilized discontinuous Galerkin (DG)
method with negative norm stabilization was proposed, which allows to
retrieve optimality in $h$ under quite weak conditions on the mesh
(allowing for the presence of very small edges).
The method there, an
hybridized formulation of which is also presented and used in the
implementation, discretizes the primal variable with polynomials of
degree $k$, and the auxiliary variable associated with the flux with
polynomials of degree $k'\in \{ k, k-1\}$, discontinuous at the
vertexes of the elements. As such choice does not satisfy the inf-sup
condition needed for the stability of the discrete problem, a
stabilization was introduced, the form of which constitutes the main
novelty of such a method. More precisely, rather that measuring the
residual term involved in a mesh dependent norm, as usually done, the
proposed stabilization makes use of a negative norm, measuring such a
residual in the space where it naturally ``lives''. This allows to
avoid the combined use of direct and inverse inequalities, which is
the main source for the lack of optimality when mesh dependent norms
are used. The negative scalar product is realized algebraically via
the introduction of an auxiliary space of minimal dimension. The
resulting formulation is shown, both theoretically and with numerical
experiments, to yield {quasi-}optimal convergence {in~$h$} even
in the presence of very small edges. However, the analysis therein is
carried out for fixed $k$, and the constants involved in the different
bounds depending on $k$. {With the proposed stabilization,}
the method {itself} %, as it is,
lacks robustness in $k$.Cockbi

In a {\em divide and conquer} approach, in this paper, we address instead
%address, independently,
the issue of the optimality with respect to $k$.  To this aim, we %then
extend the theoretical analysis of the stabilized
method {of~\cite{SilviaDaniele}} by explicitly tracking the dependence on (or independence of) the polynomial degree, and
 we present an alternative construction of the negative norm
 stabilization, which allows to achieve quasi-optimality in~$k$, this time under a stronger shape regularity assumption on
 the mesh, see Assumption~\ref{ass:meshes} (ii) and Remark~\ref{rem:latini}.
Also in this case, the negative norm is algebraically realized via the
introduction of a suitable auxiliary space. The auxiliary
space is now constructed %realized
by suitably splitting the polygonal elements into
triangles.
On each triangle, the auxiliary space is defined as the push
  forward of a space of minimal dimension, which is constructed once
  and for all by numerically solving a set of Neumann problems on a
  sufficient fine mesh on a reference triangle.

% and pushing back, on each triangle, a copy of a space of minimal
% dimension, constructed, once and for all, on a reference triangle, by
% numerically solving a set of Neumann problem on a sufficiently fine
% mesh.
Remark that, while we focus on a particular instance of the DG method, the idea of using a natural norm for the dual space in place of a mesh dependent norm used in constructing stabilization terms can be carried out to other polytopal formulations. 

\

The paper is organized as follows. In Section~\ref{sec:2}, we
  recall the stabilized DG method {from~\cite{SilviaDaniele}.
Then, in Section~\ref{sec:3}, we define the different norms and
seminorms that we will use  {in the subsequent analysis}, as well as some of
their properties{, and } prove some inverse
inequalities on polynomial spaces on the unit interval in Section \ref{sec:4}.
{In} Section~\ref{sect:theory}, we carry out {a $k$-explicit} stability and
convergence analysis {of the method; a specific} construction of a computable
bilinear form, which yields an equivalent scalar product for
$H^{-1}(K)$ on the discrete spaces used for the discretization
{in each element $K$, is
  presented} in Section~\ref{sect:stab}. {After introducing a
  hybridization of the method in Section~\ref{sec:hybrid}, which
  lesads to an efficient implementation of the method,} we present
in Section~\ref{sec:results} some numerical result confirming the validity of the theoretical estimates.

\

In the following, we will employ the notation $A \lesssim B$
(resp. $A\gtrsim B$) to indicate that $A \leq cB$ (resp. $A \geq cB$),
with $c$ positive constant independent on the mesh size parameters
$h_K$ (the diameter of the polygon $K$), $h_e$ (the length of the edge
$e$), $k$ (the polynomial degree), and possibly depending on the shape
of the polygon $K$ only via the constant in the shape regularity
Assumption \ref{ass:meshes}. We will write $A \simeq B$ to signify
that $A \lesssim B \lesssim A$. For $f \in V$, $V$ Hilbert space, and
$F \in V'$, the notation $\langle F, f \rangle$ will stand for the
action of $F$ on $f$ (the couple of dual spaces $V$ and $V'$ may vary,
its identity will be clear from the context). Moreover, in order to
avoid too cumbersome a notation, we will simply write
  ${\sup_{v}}$ instead of ${\sup_{v: v\not=0}}$ when taking the supremum over a
variable $v\not=0$ of a quantity expressed as a fraction where $v$
appears at the denominator. %, we will omit to explicitly write down the
%condition $v \not = 0$, and write $\sup_{v}$ for $\sup_{v: v\not= 0}$.

%\blu{
%\begin{itemize}
%\item In a previous paper, The method
%  there is based on a hybridized formulation, where the primal
%  variable, as well as auxiliary variable associated with ... are
%  discretized in discontinuous polynomial spaces of the same degree,
%  and a second auxiliary variable ... . The well-posedness of this two-fold mixed
%  formulation requires, in principle, that two inf-sup conditions are
%  satisfied, posing strong constraints on the choice of the
%  approximating spaces. In order to avoid this drawback, a
%  stabilization term was introduced, leading to ... (only one
%  inf-sup). This term involves ... in a negative norm. ADVANTAGES;
%  SHORT COMMENT ON ITS CONSTUCTION. 
%The resulting formulation ... (lack of $p$-robustness). 
%
%In this paper, we present an alternative construction of the negative
%norm stabilization form, in order to achieve $p$-robustness. 
%
%\item Dire qualcosa su stabilizzazioni DG / HDG.
%\item Spiegare un pochino meglio cosa vogliamo fare qui.
%\item Dire se c'e' qualcosa qui dentro che puo' servire anche a altri scopi.
%\item Outline.
%\end{itemize}
%
%Notation $\lesssim$.
%
%Notation $\langle \cdot , \cdot\rangle$.}

\section{DG method with negative norm stabilization}\label{sec:2}

As a model problem, we consider the Poisson equation with %homogeneous
Dirichlet
boundary conditions in a  polygonal domain $\Omega\subset{\mathbb R}^2$:
\begin{equation}\label{eq:Poisson}
-\Delta \w =f\quad \text{in}\ \Omega,\qquad
{\w=g} \quad\text{on}\ \partial\Omega,
\end{equation}
with $f\in L^2(\Omega)$, $g \in H^{1/2}(\partial\Omega)$.

\subsection{Assumptions on the meshes}

We consider a family $\{\Tess\}_h$ of meshes
$\Tess=\{K\}$ of the domain $\Omega$, each one
containing a finite number of polygonal elements. 
The parameter $h$ is defined as 
$h=\max_{K\in\Tess}h_K$, $h_K$ being the diameter of the polygonal
element $K$.
We denote by $\Squel$ the mesh skeleton, which is defined by
$\Squel=\bigcup_{K\in\Tess}\partial K$. For all edge $e \subset \Squel$, we let $h_e = |e|$ denote its length.

We assume that the family of meshes $\{\Tess\}_h$ satisfies the
  following properties:
\begin{assumption}\label{ass:meshes}
     There exists constants $\gamma_0,\gamma_1 > 0$ such that for all
meshes $\Tess$:
     \begin{enumerate}[(i)]
         \item  each element $K \in \Tess$ is star-shaped with respect to
a ball of radius $\geq \gamma_0 h_K$;
         \item  for each element $K$  in $\Tess$, the distance between any
two vertices of $K$ is $\geq \gamma_1 h_K$.
     \end{enumerate}
\end{assumption}

Notice that (i) and (ii) imply that there exists a constant $N>0$ such that, for all $K \in \Tess$,
 the number of edges of  $K$ is $\le N$. Moreover, it is not difficult to realize that Assumption \ref{ass:meshes} implies that $\Tess$ is graded, that is, that for all $K$, $K'$ sharing an edge it holds that $h_K \simeq h_{K'}$. 
% \st{Mainly for the sake of notational simplicity, we will from now on make a stronger assumption, namely, that the mesh is quasi uniform, with $h_K \simeq h$ for all $h$. It is however not difficult to check that all the arguments employed in the present paper carry over to the case of $\Tess$ being only graded.}
 We point out that, with Assumption~\ref{ass:meshes}, (ii), we are
   making a stronger shape regularity assumption than in~\cite{SilviaDaniele}.
 
\

Letting $\| \cdot \|_{0,D}$, $D \subset \mathbb{R}^d$, $d=1,2$, denote the standard $L^2(D)$ norm, we recall the following trace and Poincar\'e inequalities \blu{(see, e.g., \cite{beirao_stab,BrennerVEM2})}.

\subsection*{Trace inequality} Under Assumption \ref{ass:meshes}, (i),
for all $u \in H^1(K)$, we have
\begin{equation}\label{trace}
\| u \|_{0,\partial K}^2 \leq C_{tr}^2 \| u \|_{0,K} ( h_K^{-1} \| u \|_{0,K} + \| \nabla u \|_{0,K}).
\end{equation}

\subsection*{Poincar\'e inequality: 1st version} Under Assumption
\ref{ass:meshes}, (i), for all $u \in H^1(K)$, we have
\begin{equation}\label{poincare}
\inf_{q \in \mathbb{R}} \| u - q \|_{0,K} \leq C_p h_K \| \nabla u \|_{0,K}.
\end{equation}

The following result is a straightforward consequence of \eqref{poincare}.
\begin{corollary} Let $\bar u^K= |K|^{-1} \int_K u$ denote the average
  of $u$ on $K$; then, under Assumption~\ref{ass:meshes}, (i), for all $u\in
    H^1(K)$, we have
	\begin{equation}\label{corpoincare}
	 \| u - \bar u^K \|_{0,K} \leq C_p h_K \| \nabla u \|_{0,K}.
	\end{equation}
\end{corollary}

\subsection*{Poincar\'e inequality: 2nd version} Under Assumption
\ref{ass:meshes}, (i), for all $u \in H^1(K)$, setting $\widetilde u =
| \partial K |^{-1} \int_{\partial K} u$, we have
	\begin{equation}\label{poincarev2}
	\| u - \widetilde u \|_{0,K} \leq \widetilde C_p h_K \| \nabla u \|_{0,K}.
	\end{equation}
	
%	\
%	
%	Applying the trace inequality \eqref{trace} and the previous estimate to the harmonic lifting $u^\phi$ of $\phi \in H^{1/2}(\partial K)$, we obtain the following Poincar\'e-type inequality on the boundary, with $\widetilde \phi$ denoting the average of $\phi$ on $\partial K$:
%	\begin{equation}\label{poincareb}
%	\| \phi - \widetilde \phi \|_{0,\partial K} \leq \widecheck C_p h_K^{1/2} | \phi |_{1/2,\partial K}\ip{,}
%	\end{equation}
%\ip{where $| \cdot |_{1/2,\partial K}$ denote the Sobolev-Slobodeckij
%  seminorm; see Section~\ref{sect:normsedges} below. {\bf [Ho fatto un
%    search ma non ho trovato dove si usa.]}}

	% \begin{proof}
	% 	Let $\mathring u = u - \bar u$. We have
	% 	\[
	% 	u - \tilde u = \mathring u + \bar u - |\partial K |^{-1} \int_{\partial K} \mathring u - |\partial K |^{-1} \int_{\partial K} \bar u = \mathring u - |\partial K |^{-1} \int_{\partial K} \mathring u.
	% 	\]
	% 	Now we have
	% 	\begin{gather*}
	% |	\int_{\partial K} \mathring{u} | \leq | \partial K |^{1/2} \| \mathring u \|_{0,\partial K}  \leq   | \partial K |^{1/2} C_{tr} \| \mathring u \|^{1/2}_{0,K} ( h^{-1} \| \mathring u \|_{0,K} + | u |_{1,K})^{1/2} \\
	% 	\leq  | \partial K |^{1/2} C_{tr} C_p^{1/2} h^{1/2} (1+C_p)^{1/2}| u |_{1,K} \leq C_1 h | u |_{1,K}.
	% 	\end{gather*}
	%  Then we have
	% 	\[
	% 	\| u - \tilde u \|_{0,K} \leq \| \mathring u \|_{0,K} + | \int_{\partial K} \mathring u | \leq (C_p + C_1) h | u |_1.
	% 	\]

	% 	\end{proof}

\subsection{Continuous variational formulation on the
    mesh $\Tess$}

The DG methods we are going to introduce are based on the standard formulation of the primal hybrid method \cite{RaviartThomas} on $\Tess$.

Define $\w^K=\w_{|_K}$. Multiplying the equation in~\eqref{eq:Poisson} by discontinuous test
functions $v\in\prod_{K\in\Tess}H^1(K)$ 
(with an abuse of notation, we denote by $v$ also the function in $L^2(\Omega)$
such that $v_{|_K}=:v^K\in H^1(K)$)
and integrating by parts
elementwise give
\begin{equation}\label{eq:var}
\sum_{K\in\Tess}\int_K\nabla \w^K\cdot\nabla v^K -\sum_{K\in\Tess}\int_{\partial K}\nabla\w^K\cdot\mathbf{n}_K\, v^K=\int_\Omega fv,
\end{equation}
where $\mathbf{n}_K$ denotes the outer unit normal to $\partial K$. 

We define the following spaces on $\Tess$:
\[
V=\prod_{K\in\Tess}H^1(K),\qquad
\Lambda=L^2(\Squel).
\]

%\COMMENT{Qui come definire $H^{-1/2}(\Squel)$ non mi \`e chiaro. Scelte possibili sono: 1) la chiusura di $L^2$ rispetto alla norma di $\prod H^{-1/2} (\partial K)$,  2) il duale di $H^1_0(\Omega)|_{\Squel}$, 3) il sottospazio di $\prod H^{-1/2} (\partial K)$ di funzioni single valued. Altrimenti, possiamo mettere $L^2$, visto che nel nostro caso abbiamo sufficiente regolarità per farlo.}

On $\Squel$, we choose a unit normal $\normal$, taking care that, on $\partial \Omega$, $\normal$ points outwards.	
	Introduce $\la \in\Lambda$
defined as $\la=\nabla
\w_{|_{\Squel}} \cdot \normal$.
The variational formulation~\eqref{eq:var} becomes: 
find $\w\in V$, $\la \in\Lambda$
%find $u\in
%\prod_{K\in\Tess}H^1(K)$, $\lambda \in \prod_{K\in\Tess}H^{-1/2}(\partial
%K)$, and $\phi\in \{w_{|_{\Squel}}\!\!\!: w\in  H^1_0(\Omega)\}$ 
such that
\begin{equation}\label{eq:PbCont}
\begin{split}
\sum_{K\in\Tess}\int_K\nabla \w^K\cdot\nabla v^K-\sum_{K\in\Tess}\int_{\partial K}\la(\normalK\cdot \normal)v^K&=\int_\Omega
fv\qquad \forall v \in V,\\ %\in\prod_{K\in\Tess}H^1(K),\\
\sum_{K\in\Tess}\int_{\partial K}
\w^K \mu (\normalK\cdot \normal) &= \int_{\partial \Omega} g \mu\qquad
\forall \mu\in\Lambda.
\end{split}
\end{equation}
Notice that the second equation imposes the continuity of $\w$ across $\Squel$, as well as the Dirichlet boundary condition on $\partial\Omega$.

%We denote by $D: H^1(K) \to (H^1(K))'$ the operator defined as
%	\[
%	\langle Du , v \rangle = \int_K \nabla u \cdot \nabla v, \quad\text{ for all }v\in H^1(K).
%	\]
%
%Moreover, for all $K\in\Tess$, we denote by $\gamma_K^*:H^{-1/2}(\partial K)\to
%(H^1(K))^\prime$ the adjoint of the trace operator $\gamma_K:
%H^1(K)\to H^{1/2}(\partial K)$.

Observe that the well-posedness of problem~\eqref{eq:PbCont}
relies on the validity of the following  inf-sup condition:
\begin{equation*}
\inf_{\mu\in\Lambda} \sup_{v\in V}\frac{\sum_{K\in\Tess}\int_{\partial K}
v^K\mu^K}{\| \mu\|_{\Lambda}\|v\|_{V}}\ge\beta,
\end{equation*}
with a positive constant $\beta$ that, for a  suitable choice of the norms $\| \cdot \|_{\Lambda}$, $\| \cdot \|_{V}$, and under Assumption \ref{ass:meshes}, we can show to be independent of $\Tess$. 
As the well posedness of a corresponding discrete problem relies on the validity of an analogous inf-sup condition for the discrete spaces, 
a direct discretization of problem~\eqref{eq:PbCont} would require excessively 
strong assumptions on the latter. Therefore, we write a
{\em stabilized} version of problem~\eqref{eq:PbCont}.   We denote by $D: H^1(K) \to (H^1(K))'$ the operator defined as
\[
\langle Du , v \rangle = \int_K \nabla u \cdot \nabla v \quad\text{ for all }v\in H^1(K).
\]

Moreover, for all $K\in\Tess$, we denote by $\gamma_K^*:H^{-1/2}(\partial K)\to
(H^1(K))^\prime$ 
the adjoint of the trace operator $\gamma_K:
H^1(K)\to H^{1/2}(\partial K)$ and, by abuse of notation, also the operator $\gamma_K^* : \Lambda \to (H^1(K))'$ defined as
\[
\langle  \gamma_K^* \lambda, v^K \rangle = \int_{\partial K} \lambda (\normalK \cdot \normal) v^K \qquad \forall v^K \in H^1(\Omega).
\]
The abuse of notation is justified by the fact that, for $u \in
H^2(\Omega)$, if we let $\lambda \in L^2(\Squel)$ be defined as
$\lambda = \nabla u_{|_{\Squel}} \cdot \normal$, then, with the above
definition, $\gamma_K^* \lambda$ satisfies $\langle \gamma_K^*
\lambda, v \rangle = \int_{\partial K} (\partial u / \partial
\normalK) \gamma_K v$ for all $v \in H^1(K)$.
We define the jump $\Lbrack u \Rbrack$ of $u \in V$  by setting, for every interior edge $e$ shared by two elements $K^+$ and $K^-$,
\[
\Lbrack u \Rbrack = u^{K^+} \normal_{K^+} + u^{K^-} \normal_{K^-},
\]
while for $e \subset \partial \Omega \cap \partial K$ we set
\[
\Lbrack u \Rbrack = u \normalK.
\]
We observe that, for all $u \in V$, $\lambda \in \Lambda$, we have the identity
\[
\sum_{K\in\Tess} \int_{\partial K} \lambda(\normalK \cdot \normal) u^K = \int_{\Squel } \lambda \Lbrack u \Rbrack \cdot \normal
\]

\

 We consider the following stabilized problem.

\begin{problem}\label{PbContStab} Find $\w = (\w^K)_{K\in\Tess} \in
  V$, $\la \in \Lambda$ such that, for all $v = (v^K) _{K\in\Tess} \in
  V$, $\mu  \in \Lambda$,  we have
	\begin{gather}
\label{pb:contLoc1}\sum_{K\in \Tess}	\int_K \nabla \w^K \cdot\nabla v^K - \int_{\Squel } \la \Lbrack v \Rbrack \cdot \normal + t\alpha \sum_{K\in \Tess}  \scalminusone{D \w^K - \gamma_K^* \la}{D v^K} = \int_\Omega f v + t\alpha \sum_{K\in\Tess} \scalminusone{f}{D v^K}
\\
\label{pb:contLoc2}	 \int_{\Squel } \mu \Lbrack \w \Rbrack \cdot \normal - \alpha \scalminusone{Dw^K - \gamma_K^*
        \la}{\gamma_K^* \mu} = \int_{\partial \Omega} g \mu
%t
-\alpha \sum_{K\in \Tess}\scalminusone{f} {\gamma_K^* \mu},
		\end{gather}
where $\alpha$ is a positive
constant, $t=\pm 1$, and the bilinear form $\scalminusone{\cdot}{\cdot}$ denotes
the inner product in $H^1(K)'$.
	\end{problem}
It {is not} difficult to check that the same arguments we will use
in Section \ref{sect:theory} to analyze the discrete problem
{actually also} allows {us} to prove the well posedness of Problem \ref{PbContStab}.

\subsection{Discontinuous Galerkin discretization}

We define the discrete spaces
\[
V_h = \prod_{K\in\Tess} \mathbb{P}_k(K)\subset V, \qquad
\Lambda_h= \{ \lambda \in L^2(\Squel): \lambda|_e \in  \mathbb{P}_{k'}(e)\ \forall e \subset \Squel\},\]
	where $k'\in \{k-1,k\}$,
	 and where $\mathbb{P}_k(K)$ 
	 (resp. $\mathbb{P}_{k'}(e)$)
	 denotes the space of
polynomials  of degree at most $k$ 
	(resp. $k'$)
	 in two variables restricted to $K$ 
	(resp. $e$).

%\COMMENT{Is it worth remarking that the element of $\Phi_h$ are univariate piecewise polynomials in the curvilinear abscissa, since the restriction of a bivariate polynomial to a straight line is a univariate polynomial? Or maybe we can observe it when we prove the inverse inequalities}

The discrete version of the stabilized problem~\eqref{PbContStab}
reads as follows.

\begin{problem}\label{PbGlob} Find $\uh = (\uh^K)_{K\in\Tess} \in
  V_h$, $\lah  \in \Lambda_h$ such that, for all $v = (v^K)
  _{K\in\Tess} \in V_h$ and $\mu  \in \Lambda_h$, we have
	\begin{multline}
		\label{eq:PbGlob-loc1}
		\sum_{K\in\Tess}	\int_K \nabla \uh^K \cdot\nabla v^K - \int_{\Squel } \lah \Lbrack v \Rbrack \cdot \normal + t\alpha \sum_{K\in\Tess} 
s_K (D \uh^K - \gamma_K^* \lah,D v^K)
\\ = \int_\Omega f v + t\alpha \sum_{K\in\Tess} s_K (f,D v^K)\end{multline}	
		\begin{gather}
		\label{eq:PbGlob-loc2}	 \int_{\Squel } \mu \Lbrack \uh \Rbrack \cdot \normal - \alpha \sum_{K\in\Tess} s_K (D\uh^K - \gamma_K^*
		\lah, \gamma_K^* \mu) = \int_{\partial \Omega} g \mu
		%t
		-\alpha \sum_{K\in\Tess}s_K (f, \gamma_K^* \mu).
		\end{gather}
Here, $s_K: (H^1(K))' \times (H^1(K))' \to \mathbb{R}$ is a continuous 
bilinear form that, when restricted to %average free
elements $F \in \gamma_K^*(\Lambda_h)$ with $\langle F, 1 \rangle = 0$,  is spectrally
equivalent to the $(H^1(K))'$ inner product $\scalminusone\cdot\cdot$.
	\end{problem}

%\subsection{Option 1. Continuous $\phi$

More precisely,
we define by duality the following norm and seminorm for elements
$F\in (H^1(K))^\prime$ (see Section \ref{sec:3} for more details):
\begin{equation}\label{dualKseminorm}
\| F \|_{-1,K} = \sup_{g\in H^1(K)}\frac{\langle
	F,g\rangle}{(| K |^{-2}| \int_K g |^2 + \| \nabla g\|^2_{0,K})^{1/2}}, \quad
|F|_{-1,K}=\sup_{g\in H^1(K),\ \int_K g=0}\frac{\langle
  F,g\rangle}{\| \nabla g\|_{0,K}},
\end{equation}
and we make the following assumptions on the stabilization forms $s_K (\cdot,\cdot)$.

\begin{assumption}\label{sK1} (Continuity) There exists a constant $M(k)>0$,
  possibly depending on $k$, such that
	\[
	s_K(F,G) \leq M(k) | F |_{-1,K} | G |_{-1,K}\qquad\forall F,G
        \in (H^1(K))', \ \forall K\in\Tess.
	\]
\end{assumption}
\begin{assumption}\label{sK2} (Coercivity) There exists a constant $\rho(k)>0$,
  possibly depending on $k$, such that
	\[
	s_K(\gamma_K^* \lambda,\gamma_K^* \lambda ) \geq   \rho(k) |
        \gamma_K^* \lambda |_{-1,K}^2
\qquad\forall \lambda \in \Lambda_h\ ,\forall K\in\Tess.
	\]	
\end{assumption}
Notice that $\rho(k) \leq M(k)$.

\

Computable stabilization forms $s_K (\cdot,\cdot): (H^1(K))' \times (H^1(K))'\to {\mathbb R}$ satisfying
Assumptions~\ref{sK1} and~\ref{sK2} will be introduced below
in Section~\ref{sect:stab}. Before presenting a stability and
error analysis of the DG formulation in Problem~\ref{PbGlob} (see
Section~\ref{sect:theory} below), we introduce some norms and
seminorms, together with their
properties (Section~\ref{sec:3}), and recall some properties of
polynomial spaces (Section~\ref{sec:4}).
%as well as some relevant inequalities. 

\section{Norms and seminorms}\label{sec:3}

We start by defining local norms and seminorms on $e\subset\Squel$, on
$K$ and  $\partial K$, $K\in\Tess$, that are convenient in
the application of scaling arguments, in particular when negative norms are concerned,
and define global norms on $\Tess$ and $\Squel$.

\subsection{Norms and seminorms on elements}\label{sect:normselements}

We define the following norms and seminorms for the Sobolev spaces
$H^1(K)$ and its dual $(H^1(K))^\prime$.

For $u\in H^1(K)$, we let
\begin{equation}\label{defKnorm}
\| u \|^2_{1,K} = | \bar{u}^{K} |^2 + | u |^2_{1,K},
\end{equation}
with
\[
\bar{u}^{K} = \frac{1}{|K|} \int_K u, \qquad | u |_{1,K}^2 = \int_K | \nabla u |^2.
\]
% For $F\in (H^1(K))^\prime$ we define
% \begin{equation}\label{dualKseminorm}
% |F|_{-1,K}=\sup_{g\in H^1(K),\ \int_K g=0}\frac{\langle F,g\rangle}{|g|_{1,K}}.
% \end{equation}

%Moreover we define a norm and a seminorm for $H^1(K)$
%\[
%\| F \|_{-1,K} =\sup_{u \in H^1(K)} \frac{\langle F, u \rangle}{\| u \|_{1,K}}, \qquad 
%| F |_{-1,K} = \sup_{u \in H^1(K):\ \int_K u = 0} \frac{\langle F, u \rangle}{| u |_{1,K}},
%\]

Observe that, if $\| \cdot \|_{-1,K}$ and $| \cdot |_{-1,K}$  are defined by duality with the norm $\| \cdot \|_{1,K}$ as in \eqref{dualKseminorm},  then it holds (see \cite{SilviaDaniele})
\[
\| F \|_{-1,K}^2  = | \langle F,1 \rangle|^2 + | F |_{-1,K}^2. 
\] 

\

Under Assumption \eqref{ass:meshes}, we have the following proposition.

\newcommand{\pezzobordo}{\tau}

	\begin{proposition}\label{prop:genPW} Let $\pezzobordo \subseteq \partial K$, $K \in \Tess$, be a connected subset of $\partial K$ with $| \tau | \geq \gamma_1 h_K$. Then,  for all $u \in H^1(K)$,
		letting $\bar u^\pezzobordo = | \pezzobordo|^{-1} \int_{ \pezzobordo} u$, we have
	\[
	\| u \|_{1,K}^2 \simeq  | \bar u^\pezzobordo  |^2 + | u |^2_{1,K}.
	\]	
\end{proposition}

\newcommand{\Cpoinctau}{\widehat C_{p}}

	\begin{proof}  We have, with $\bar u^K = |K|^{-1} \int_K u$,
		\begin{multline}
		\int_K | u(x) |^2 \,dx = \int_K | u(x) - \frac 1 {|\pezzobordo|} \int_{ \pezzobordo} u(s)\, ds |^2 \,dx \leq \int_K \frac 1{| \pezzobordo|} \int_{ \pezzobordo} | u(x) - u(s) |^2\,ds \,dx\\
		\leq  \int_K \frac 2{| \pezzobordo|} \int_{ \pezzobordo} | u(x) - \bar u^K|^2\,ds \,dx +  \int_K \frac 2{| \pezzobordo|} \int_{ \pezzobordo} | u(s) - \bar u^K |^2\,ds \,dx \\
		\leq
		2 \| u - \bar u^K \|_{0,K}^2 + 2 \frac {2|K|} {\pezzobordo}\| u - \bar u^K \|^2_{0,\pezzobordo} \leq 	2 \| u - \bar u^K \|_{0,K}^2 + 2 \frac {|K|} {\pezzobordo}\| u - \bar u^K \|^2_{0,\partial K}.
		\end{multline}
		Then, by applying the trace inequality~\eqref{trace}
                and the Poincar\'e inequality \eqref{poincare}, for
                all $u\in H^1(K)$ with $\int_{\pezzobordo} u = 0$, we have
		\begin{equation}\label{eq:PW3}
		\| u \|_{0,K} \leq \Cpoinctau h_K | u |_{1,K}
		\end{equation}
		with $\Cpoinctau$ only depending on the shape
                regularity constants $\gamma_0$ and $\gamma_1$. 
		Conversely, by applying once again the trace inequality~\eqref{trace} and the Poincar\'e inequality \eqref{poincare}, we can write 
		\[
		| \bar u^\pezzobordo |^2 \leq | \pezzobordo|^{-1} \| u \|^2_{0,\pezzobordo} \leq h_K^{-1} \| u \|^2_{0,\partial K} 
		\leq h_K^{-2} \| u - \bar u^K  \|_{0,K}^2 + | \bar u^K |^2 + | u |^2_{1,K} \lesssim \| u \|_{1,K}^2.
		\]
\end{proof}

\subsection{Norms and seminorms on edges}\label{sect:normsedges}
%Let us assume that the edge $e$ has length $h$.

For every edge $e$ of the mesh, we define the following norms and seminorms for the Sobolev spaces
$H^{s}(e)$ and their dual spaces~$(H^{s}(e))^\prime$,
  $0<s<1$. 
% $H^{-1/2}(e)$. 

For $\phi \in H^{s}(e)$, we let 
\begin{equation}\label{defedgenorm}
\| \phi \|_{s,e}^2 = |e|^{1-2s} | \bar \phi |^2 + | \phi |^2_{s,e},
\end{equation}
with
\[
\bar \phi = \frac 1 {|e|} \int_e \phi,  \qquad | \phi |_{s,e} = \int_e \int_e \frac{| \phi(x)-\phi(y)|^2}{| x - y |^{2s+1}}\,dx\,dy.
\]
On $(H^{s}(e))^\prime$, we define %the following norm and a seminorm: 
\begin{equation}\label{defedgedualnorm}
\| \lambda \|^2_{-s,e} = |e|^{2s-1} | \langle \lambda, 1 \rangle |^2 + | \lambda
|^2_{-s,e}, \qquad | \lambda |_{-s,e} = \sup_{{\phi \in H^{s}(e)}:{\int_e \phi = 0}}\, \frac{ \langle \lambda , \phi \rangle }{| \phi |_{s,e}}.
\end{equation}

\newcommand{\uI}{I}

%	\COMMENT{
%Other scale relations that may, or may not, be used in the following are, for $0 \leq s \leq 1$
%\[
%| \phi |_{s,e} = | e|^{1/2 - s} | \widehat \phi |_{s,\uI}, \qquad | \phi |_{-s,e} = | e |^{1/2+s} | \widehat \phi |_{-s,\uI}.
%\]	
%}

\

The two norms defined by (\ref{defedgenorm}) and
(\ref{defedgedualnorm}), respectively, satisfy the duality relations
	\begin{gather*}
	\| \lambda \|_{-s,e} = \sup_{\phi\in H^{s}(e)} \frac { \int_e \lambda \phi }{\| \phi \|_{s,e}},\qquad
\| \phi \|_{s,e} = \sup_{\lambda\in (H^{s}(e))'} \frac { \int_e \lambda \phi }{\| \lambda \|_{-s,e}},
\end{gather*}
see~\cite[Lemma~2.1]{SilviaDaniele}.

\

\newcommand{\vphi}{\varphi}
\newcommand{\veps}{\varepsilon}

\blu{

	On $e:=(a,b)$, we will also consider the spaces
	$H^s_0(e)$ ($s \in ]0,1[$, $s \not= 1/2$) and
	$H^{1/2}_{00}(e)$ of functions whose extension by zero is in
	$H^s(\mathbb{R})$ ($s \not= 1/2$) and
	$H^{1/2}(\mathbb{R})$ respectively, which we will equip with the norms
	\begin{gather}
	\| \vphi \|^2_{H^{s}_{0}(e)} = | \vphi |^2_{H^{s}(e)} +
	\int_e \frac {| \vphi(x) |^2}{| x - a |^{2s}} \,ds(x) +\int_e \frac
	{| \vphi(x) |^2}{| x - b |^{2s}} \,ds(x),\qquad s\not =
        1/2, %0
        \\
		\| \vphi \|^2_{H^{1/2}_{00}(e)} = | \vphi |^2_{H^{1/2}(e)} +
		\int_e \frac {| \vphi(x) |^2}{| x - a |} \,ds(x) +\int_e \frac
		{| \vphi(x) |^2}{| x - b |} \,ds(x).
	\end{gather}
	For $s = 1$, we set $\| \vphi \|_{H^1_0(e)} = | \vphi |_{1,e}$.

	\
	
	We recall that, for $s < 1/2$, the two spaces $H^s(e)$ and $H^s_0(e)$
	coincide, and
	the two corresponding norms are equivalent. However, the constant in the
	equivalence depends on $s$ and it
	explodes as
	$s$ converges to $1/2$. The behavior of such constant as $s$ approaches the
	limit value $1/2$ is given by the following bound, which holds for all
	$\zeta \in H^ {1/2-\veps}(e)$, with $\bar \zeta = | e |^{-1} \int_e \zeta$ (see \cite{BertoluzzaMathComp}):	
\begin{equation}
\| \zeta \|_{H^{1/2-\veps}_0(e)} \lesssim
\frac 1 \veps
| \zeta |_{1/2-\veps, e} + \frac {|e|^\veps} {\sqrt{\veps}}
|\bar \zeta |.\label{penultima}
\end{equation}

By a simple duality argument, it is not difficult to check that we have
		\begin{equation}
\| \zeta \|_{-1/2+\veps, e} 	 \lesssim
\frac 1 \veps \| \zeta \|_{(H^{1/2-\veps}_0(e))'}.\label{penultimadual}
\end{equation}

}

Observe that the seminorm $| \cdot |_{1/2,e}$  and the norm $\| \cdot \|_{H^{1/2}_{00}(e)}$ are scale invariant. In
fact, letting $\uI = [0,1]$ and $e = [0, |e|]$, and setting $\widehat x = |e|^{-1}x$,  for $\widehat \phi \in H^{1/2}(\uI)$ (resp. $\widehat \phi\in H^{1/2}_{00}(\uI)$) and $\phi (x) = \widehat \phi(\widehat x) \in H^{1/2}(e)$  (resp. $\phi (x) = \widehat \phi(\widehat x)\in H^{1/2}_{00}(e)$), we have the identity
\begin{equation}\label{scalinghalfnorm}
| \widehat \phi |_{1/2,\uI} = | \phi |_{1/2,e} \qquad \text{(resp.}\ \|\widehat \phi\|_{H^{1/2}_{00}(\uI)} = \| \phi \|_{H^{1/2}_{00}(e)}  \text{)}.
\end{equation}

\

For the $| \cdot |_{-1/2,e}$ seminorm and the $\| \cdot \|_{(H^{1/2}_{00}(e))'}$ norm, for $\widehat \lambda \in
L^2(\uI)$ and $\lambda (x) = \widehat \lambda(\widehat x)$, we instead have 
\[
| \lambda |_{-1/2,e} = |e| | \widehat \lambda |_{-1/2,\uI}, \qquad 
\| \lambda \|_{(H^{1/2}_{00}(e))'} = | e | 
\| \widehat \lambda \|_{(H^{1/2}_{00}(\uI))'}.
\]
In fact, 
\begin{gather*}
| \lambda |_{-1/2,e} = \sup_{{\phi \in H^{1/2}(e)}:\ {\int_e \phi = 0}} \frac{ \int_e \lambda(x) \phi(x) \,dx}{| \phi |_{1/2,e}} = |e|
\sup_{{\widehat \phi \in H^{1/2}(\uI)}:\ {\int_{\uI} \widehat \phi= 0}} \frac{ \int_{\uI} \widehat \lambda(\widehat x) \widehat \phi(\widehat x) \, d\widehat x }{| \widehat \phi |_{1/2,\uI}},
\end{gather*}
(and analogously for the $(H^{1/2}_{00}(e))'$ norm).

The norm $\| \cdot \|_{H^{1/2}_{00}(e)}$ controls the norm $\|\cdot \|_{1/2,e}$ uniformly in $|e|$, namely
\[
\| \phi \|_{1/2,e} \lesssim \| \phi  \|_{H^{1/2}_{00}(e)}.
\]
This readily follows from
\[
| \bar \phi |^{2} \leq \frac{1}{|e|} %|e| 
\int_e | \phi (x) |^{2} \,dx \leq  \int_e\frac{| \phi (x) |^{2} }{|x|} \,dx.
\]

By duality, %(recall that $H^{-1/2}(e) = (H^{1/2}_{00}(e))'$), 
we have that
\begin{equation}\label{compdualuno}
\| \lambda \|_{(H^{1/2}_{00}(e))'} \lesssim \| \lambda \|_{-1/2,e}.
\end{equation}

\subsection{Norms and seminorms on element boundaries}\label{sect:normsboundaries}

We define the norm in $H^{1/2}(\partial K)$ as
\[
\|\phi\|_{1/2,\partial K}^2=|\bar{\phi}^{\partial K}|^2+|\phi|_{1/2,\partial K}^2,
\]
with
\[
\bar{\phi}^{\partial K} =\frac{1}{| \partial K|} \int_{\partial K} \phi,\qquad
| \phi |_{1/2,\partial K} = \int_{\partial K} \int_{\partial K} \frac{| \phi(x)-\phi(y)|^2}{| x - y |^2}\,dx\,dy.
\]

\

%With this norm
We have the following equivalence between this norm %the $H^{1/2}(\partial K)$ 
and the norm obtained via the trace operator: 
	\begin{equation}\label{tracenorm}
	\inf_{u \in H^1(K):\ u=\phi \text{ on }\partial K} \| u \|_{1,K} \simeq \| \phi \|_{1/2,\partial K}.
	\end{equation}
In fact, by using the definition of $\| \cdot \|_{1/2,\partial
    K}$ and, recalling that
\begin{equation}\label{traceseminorm}
|\phi|_{1/2,\partial K} \simeq \inf_{u \in H^1(K):\ u=\phi \text{ on
  }\partial K} | u |_{1,K},
\end{equation} we can write
	\[
	\| \phi \|_{1/2,\partial K}^2 = | \bar{\phi}^{\partial K} |^2
        + | \phi |_{1/2,\partial K }^2 \simeq 
\inf_{u \in H^1(K):\ u=\phi \text{ on
  }\partial K}(| \bar u^{\partial K} |^2 + | u |^2_{1,K }),
	\]
and Proposition \ref{prop:genPW} implies the equivalence.

\

We now state some relations between dual norms on
$\partial K$. %\blu{Skip proof and quote the other paper?}

\begin{proposition}\label{prop:lambdalambda}
Let $\lambda \in L^2(\Squel)$. Under Assumption
\ref{ass:meshes}, for all $K \in \Tess$, we have
	\[
	\left(\sum_{e\subset\partial K} \| \lambda\|^2_{(H^{1/2}_{00}(e))'}	\right)^{1/2} \lesssim \| \gamma_K^* \lambda \|_{-1,K} \lesssim \left(\sum_{e\subset\partial K}\| \lambda \|^2_{-1/2,e}\right)^{1/2},
	\]
\end{proposition}

\begin{proof} We start by proving the second bound. Let $\lambda \in L^2(\partial K)$. We have 
	\begin{gather*}
	\sup_{u \in H^1(K)} \frac{\int_{\partial K} \lambda u }{\| u
          \|_{1,K}} = 	\sup_{u \in H^1(K)}
        \frac{\sum_{e\subset\partial K}\int_{e} \lambda u }{\| u
          \|_{1,K}} 
%\leq 	\sup_{u \in    H^1(K)}\sum_{e\subset\partial K} \frac{\int_{e} \lambda u }{\| u \|_{1,K}} 
\leq 	\sup_{u \in H^1(K)}\sum_{e\subset\partial K} \frac{\| \lambda \|_{-1/2,e} \| u \|_{1/2,e}}{\| u \|_{1,K}}.
	\end{gather*}
	Let us now compare $\| u \|_{1/2,e}$ with $\| u
        \|_{1,K}$. From the definition of $\| u \|_{1/2,e}$ and Proposition~\ref{prop:genPW}, using \eqref{traceseminorm} we have, with $\bar u^e = | e |^{-1} \int_e u$,
        \begin{equation*}
        	\| u \|^2_{1/2,e} = | \bar u^e |^2 + | u |^2_{1/2,e} \lesssim  | \bar u^e |^2 + | u |_{1/2,\partial K} 
        	\lesssim
        	  | \bar u^{e} |^2 + | u |^2_{1, K}, \lesssim \| u \|_{1,K},
        \end{equation*}
	yielding
	\[ \sup_{u \in H^1(K)} \frac{\int_{\partial K} \lambda u }{\|
            u \|_{1,K}} \lesssim \sum_{e\subset\partial K}\| \lambda
          \|_{-1/2,e}
          \lesssim \left(
          \sum_{e\subset\partial K}\| \lambda
          \|^2_{-1/2,e}
          \right)^{1/2},
          \]
where we used that the number of edges of $K$ is uniformly bounded
thanks to Assumption~\ref{ass:meshes}. This proves the second
  bound of the statement.

	As far as the first bound is concerned, we remark that \eqref{tracenorm} implies that
	\[
	\| \phi \|_{1/2,\partial K}^{-1} \lesssim \sup_{u \in H^1(K): \ u|_{\partial K}=\phi} \| u \|_{1,K}^{-1}.
	\]
	
% we start by
%         introducing the space $\Phi_0 \subset H^{1/2}(\partial K)$ defined as
% 	\[
% 	\Phi^0 = \{\phi \in H^{1/2}(\partial K): \phi|_{e} \in
%         H^{1/2}_{00}(e) \ \forall e\subset\partial K\},
% 	\]
% 	and we observe that for all $\phi \in \Phi_0$ we have
% 	\[
% 	\Big| \frac{1}{|e|} \int_e \phi \Big| \lesssim \| \phi \|_{H^{1/2}_{00}(e)}
% 	\]
% 	(prove it on the reference element).
% 	Moreover we have
% 	\[
% 	|\partial K|^{-1} | \int_{\partial K}\phi | \leq \max_e \frac{|e|}{|\partial K|} \sum_e |e|^{-1} | \int_e \phi | \lesssim  \max_e \frac{|e|}{|\partial K|} \sum_e \| \phi \|_{H^{1/2}_{00}(e)},
% 	\]
% 	and
% 	\[
% 	| \phi |_{H^{1/2}(\partial K)} \leq \sum_{e} \| \phi
%         \|_{H^{1/2}_{00}(e)}.
% 	\]
% Combining these two inequalities, we obtain
% \[
% \| \phi \|_{H^{1/2}(\partial K)} \le \sum_{e} \| \phi
%         \|_{H^{1/2}_{00}(e)}.
% \]
% 	Therefore, 
%
Then we can write
	\begin{multline*}
\| \lambda \|_{(H^{1/2}_{00}(e))'}
=\sup_{\phi\in H^{1/2}_{00}(e)}\frac{\int_{e}\lambda\phi}{\|\phi\|_{H^{1/2}_{00}(e)}}
=\sup_{\phi\in H^{1/2}(\partial K),\ \phi_{|_{\partial K\setminus e}}=0}
\frac{\int_{\partial
      K}\lambda\phi}{\|\phi\|_{1/2,\partial K}}
\le \sup_{\phi\in H^{1/2}(\partial K)}\frac{\int_{\partial
      K}\lambda\phi}{\|\phi\|_{1/2,\partial K}}\\
\lesssim \sup_{\phi \in H^{1/2}(\partial K)} \sup_{u\in H^1(K): \ u|_{\partial K} = \phi}
\frac{\int_{\partial
      K}\lambda u}{\|u\|_{1,K}} =
\sup_{u\in H^1(K)}\frac{\int_{\partial
      K}\lambda u}{\|u\|_{1,K}} = \| \gamma_K^* \lambda \|_{-1,K}.
	\end{multline*}
	By squaring and adding up the contributions of the different
        edges, taking once again into account that the number of edges is
        uniformly bounded, we obtain the first bound, and the
          proof is complete.
\end{proof}

\subsection{Global norms and seminorms}\label{sec:globalnorms}

\newcommand*{\mvl}[1]{\{\!\!\{#1\}\!\!\}}             % mean value

%We introduce the standard DG notation for the average and jump of a
%piecewise continuous function $w$ on $\Tess$ across an interior edge
%$e\subset\Squel$ shared by the two elements $K_1$ and $K_2$:
%\[
%\mvl{w}=\frac{1}{2}w^{K_1}+\frac{1}{2}w^{K_2}, \qquad
%\Lbrack w \Rbrack=w^{K_1}{\mathbf n}_{K_1}+w^{K_2}{\mathbf n}_{K_2},
%\]
%where ${\mathbf n}_K$ denotes the unit normal vector to $\partial K$ pointing
%outside $K$. 
%On $e\subset\partial \Omega$  boundary edge  belonging to the element $K$, 
%we set $\mvl{w}=w^{K}$ and $\Lbrack w \Rbrack=w^{K}{\mathbf n}_{K}$.

We define the following global seminorms and norms on $\Tess$ and $\Squel$: 

\begin{align*}
|u|_{1,\Tess}^2&=\sum_{K\in\Tess} | u^K |_{1,K}^2, &&\forall u\in \prod_{K\in\Tess}H^1(K),\hfill\\
\| u \|_{1,\Tess}^2 &= |u|_{1,\Tess}^2+ \sum_{e\subset
                  \Squel} | \Lbrack \bar u \Rbrack |^2,
&&\forall u\in \prod_{K\in\Tess}H^1(K),\hfill\\
\tb u \tb_{1,\Tess}^2 &=\sum_{K\in\Tess} \| u^K \|_{1,K}^2, && \forall u \in \prod_{K\in\Tess}H^1(K),\hfill\\
\| \lambda \|_{-1/2,\Squel}^2 &=
                                \sum_{e\subset\Squel} \| \lambda^e
                                \|^2_{-1/2,e}, &&\forall
                                \lambda\in
                                                \prod_{e\subset\Squel}(H^{1/2}(e))', \hfill\\
 \tb \lambda \tb^2_{-1/2,\Squel}, &= \sum_{K\in\Tess} \| \gamma_K^*
                                   \lambda^K \|^2_{-1,K},
&&\forall \lambda\in \Lambda, \hfill
\end{align*}
where the superscripts $K$ and $e$ denote the restrictions to $K$
  and $e$, respectively, and
$\bar u\in\prod_{K\in\Tess} H^1(K)$ is such that $\bar{u}^K=\frac{1}{|K|}\int_K u$.
%denotes the piecewise constant function defined on each
%$K\in\Tess$ as $\bar u=\frac{1}{|K|}\int_K u$. 
Notice that $|\Lbrack \bar u
\Rbrack|=|\bar{u}^{K^+}- \bar{u}^{K^-}|$ on $e$, if $e$ is an
interior edge shared by the elements $K^+$ and $K^-$, or $|\Lbrack \bar u
\Rbrack|=|\bar{u}^{K}|$, if $e$ is a boundary edge that belongs to the
element $K$.

\

For all $u\in \prod_{K\in\Tess}H^1(K)$, the following Poincar\'e-type inequality holds true~(see~\cite[Lemma~2.6]{SilviaDaniele} with $H_e = h_e = | e |$):
\[
\| u \|_{0,\Omega} \lesssim \| u \|_{1,\Tess}.
\]
Moreover, it is easy to check that
\begin{equation}\label{stimanormath1}
\| u \|^2_{1,\Tess} \lesssim | u |_{1,\Tess}^2 + h^{-2} \| u \|_{0,\Omega}^2.
\end{equation}
%	\[
%	\| u \|_{0,\Omega}^2 \lesssim \it  | u^K |_{1,K}^2 + \sum_e | \Lbrack \bar u \Rbrack |^2
%	\]

%\blu{[The lemma with the sketch of the proof is now commented.]}

% \begin{lemma}[Poincar\'e] Let $u \in \prod_K H^1(K)$ and let $\check u = (\check u^K)$ denote the piecewise constant function with $\check u^K$ denoting the average of $u^K$ on $K$. Then

% \end{lemma}

% \begin{proof} We can write
% 	[.. bla bla ..]
	
% 	Let $w$ be the solution of  
% 	\[
% 	-\Delta w = \check u, \text{ in }\Omega, \qquad w = 0, \text{ on }\partial \Omega. 
% 	\]
% 	(23) of the three fields paper holds (same proof). Observe that $w \in H^2$ implies continuity of the normal derivative across the skeleton. We can define 
% 	\[
% 	\bar \mu^e = \frac{1}{| e| }\int_e \frac {\partial w} {\partial \nu_e }. 		\]
% 	We have (24) of three fields paper, but also (same proof)
% 	\[
% 	\| \bar \mu^e \|_{-1/2,e} \leq \| \partial w/\partial \nu^e \|_{-1/2,e}
% 	\]
% 	Then we write
% 	\begin{gather*}
% 	\| \check u \|^2_{0,\Omega} = - \it  \int_{\partial K} \frac {\partial w} {\partial \nu^K} \check u^K ds =
% 	\sum_e \int_e \bar \mu^e {\Lbrack} \check u  {\Rbrack} \leq \sum_e h_e | \bar \mu^e | | \Lbrack \check u \Rbrack | \\
% 	\lesssim 
% 	(\sum_e h_e^2 | \bar \mu^e |^2)^{1/2}  (\sum_e  | \Lbrack \check u \Rbrack |^2)^{1/2} 
% 	\end{gather*}
	
% \end{proof}

%\newpage

\section{Inverse inequalities in polynomial spaces on the unit interval}\label{sec:4}

\newcommand{\Pko}{P_{k+2}^0}

In this section, we recall some inverse inequalities for polynomials in
positive Sobolev
norms, and establish inverse inequalities in negative Sobolev
norms. These results will be used in Section~\ref{sect:stab}.

Assume that $I$ is an interval of unit length. 
We start by recalling that, for all
$p \in \mathbb{P}_k(I)$, it holds that
\begin{equation}\label{inversebase}
\|p\|_{1,I}\lesssim k^2 \|p\|_{0,I},
\end{equation}
\blu{see, e.g., \cite[Theorem~3.91]{SchwabBook}}. 

%\COMMENT{ 
%Queste non si usano mai, e sono abastanza standard. Io le toglierei	
%	By a standard space
%interpolation argument, we immediately get  
%\begin{equation}\label{eq:interpmezzo}
%\|p\|_{1/2,I}\lesssim k \|p\|_{0,I}, \qquad \text{ for all } p \in \mathbb{P}_k(I).
%\end{equation}
%as well as 
%	\begin{equation}\label{eq:interpmezzo0}
%	\|p\|_{H^{1/2}_{00}(I)}\lesssim k \|p\|_{0,I}, \qquad \text{ for all } p \in \mathbb{P}_k(I) \cap H^1_0(I).
%	\end{equation}
%	Do we need this last one?}

We prove now inverse inequalities in negative Sobolev norms.

\begin{lemma}\label{lem:inversenegative}
Let $0 \leq r \leq s \leq 1$. 
If $r,s\ne 1/2$,
for all $p \in \mathbb{P}_k(I)$, we have 
\begin{equation}\label{eq:inversenegative}
\| p \|_{(H^r_0(I))'} \lesssim k^ {2(s-r)} \| p \|_{(H^s_0(I))'}.
\end{equation}
Moreover, for  $1/2 \leq s \leq 1$ and $0 \leq r \leq 1/2$, for all $p \in \mathbb{P}_k(e)$, we have 
\begin{equation}\label{eq:inversenegativemezzo}
\| p \|_{(H^{1/2}_{00}(I))'} \lesssim k^ {2s-1} \| p \|_{(H^s_0(I))'}, \qquad 
\| p \|_{(H^r_0(I))'} \lesssim k^ {1-2r} \| p \|_{(H^{1/2}_{00}(I))'}.
\end{equation}
\end{lemma}

\begin{proof}
We let $\Pko:L^2(\uI) \to \mathbb{P}_{k+2}^0(\uI) = \mathbb{P}_{k+2}(\uI) \cap H^1_0(\uI)$ be defined as
\[
\int_\uI (\Pko \phi - \phi) q = 0 \quad \text{ for all } q \in \mathbb{P}_k(\uI).
\]
It is easy to see  that $\Pko$ is well defined. Indeed, for all $p \in \mathbb{P}_{k+2}^0(I)$, letting $q =- p''\in \mathbb{P}_{k}(I)$, we have
\[
\int_\uI p q = \int_\uI | p' |^2 \not = 0.
\]
As  $\dim(\mathbb{P}_{k+2}^0(I)) = \dim(\mathbb{P}_k(I))$, this implies that $\Pko$ is well defined.
%
%We have (see [Bernardi, Maday, Patera, preprint Mortar] appendix B) for all $\phi \in H^1_0(\uI)$
%\[
%\| \Pko  \phi \|_{1,\uI} \lesssim \| \phi \|_{1,\uI}.
%\]
We can write:
\begin{equation}\label{eq:proofinv1}
\| p \|_{0,\uI}^2 = \int_\uI | p |^2 = \int_\uI p \Pko(p) \leq \| p \|_{(H^1_0(\uI))'} \| \Pko (p) \|_{H^1_0(\uI)}.
\end{equation}
We then need to bound $\| \Pko(p) \|_{H^1_0(\uI)}$. We have
\[
\| \Pko(p) \|_{H^1_0(\uI)} \lesssim  \int_\uI |  \Pko(p)' |^2 = - \int_\uI \Pko(p)'' \Pko(p)  = 
 - \int_\uI \Pko(p)'' p,
\]
where the last identity stems from the definition of the projector $\Pko$, as $\Pko(p)'' \in \mathbb{P}_k(\uI)$. Then, using \eqref{inversebase}, we have
\[
\| \Pko(p) \|_{H^1_0(\uI)} \lesssim \| \Pko(p)'' \|_{0,\uI} \| p \|_{0,\uI} \lesssim k^2 | \Pko(p) |_{1,\uI} \| p \|_{0,\uI},
\]
whence by dividing both sides by $\| \Pko (p) \|_{H^1_0(\uI)} = | \Pko (p) |_{1,\uI}$ and substituting in \eqref{eq:proofinv1}, we obtain
\[
\| p \|_{0,\uI}  \lesssim k^2    \| p \|_{(H^1_0(\uI))'}.
\]

By interpolating between $(H^1_0(I))'$ and  $L^2(I)$, we obtain that
\[
\| p \|_{0,I} \lesssim k^{2s} \| p \|_{(H^s_0(I))'}.
\]
Finally, by interpolating between $(H^s_0(I))'$ and $L^2(I)$, we
get~\eqref{eq:inversenegative}. As the space $[(H^1_0(\uI))',L^2(\uI)]_{\theta}$ obtained by space interpolation between $(H^{1}_0(\uI))'$ and $L^2(\uI)$ is, for $\theta=1/2$, $(H^{1/2}_{00}(\uI))'$ rather than $(H^{1/2}_{0}(\uI))'$, for either $s=1/2$ or $r=1/2$ the above argument gives us~\eqref{eq:inversenegativemezzo}.
\end{proof}

For polynomial functions, we also have the following lemma which, combined with the bound \eqref{compdualuno}, states the equivalence of the norms for $(H^{1/2}(e))'$ and $(H^{1/2}_{00}(e))'$
\begin{lemma}\label{lem:log} For all $p \in \mathbb{P}_k(I)$  we have
\begin{equation}\label{eq:wish1}
\|\lambda\|_{-1/2,e}\lesssim \log k\, \|\lambda\|_{(H^{1/2}_{00}(e))'}.
\end{equation}
	\end{lemma}

\begin{proof} Denoting by $\bar\lambda^{e_i}$ the average of $\lambda$ on
$e$, we have
\[
\|\lambda\|^2_{-1/2,e} = \left|\int_{e}\lambda \right|^2 + |\lambda|^2_{-1/2,e}.
\]
We bound the two terms on the right-hand side separately.
For the first one, using \eqref{eq:inversenegative} and \eqref{penultima}, we have
\begin{equation*}
\left|\int_{e}\lambda\right|\lesssim|\lambda|_{(H^{1/2-\varepsilon}_0(e))'}
|1|_{H^{1/2-\varepsilon}_0(e)}\\ \lesssim
\frac{k^{2\varepsilon}}{h_{e}^{\varepsilon}}
|\lambda|_{(H^{1/2}_{00}(e))'}\frac {h_{e}^\varepsilon} {\sqrt{\varepsilon}} =\frac{k^{2\varepsilon}}{\sqrt{\varepsilon}}|\lambda|_{(H^{1/2}_{00}(e))'}.
\end{equation*}

%\COMMENT{Check scaling of norms.}

For the second term, we have
\begin{gather*}
|\lambda|_{-1/2,e}=\sup_{\phi\in
	H^{1/2}(e_i)}\frac{\int_{e}\lambda\phi}{|\phi|_{1/2,e}}
\le\sup_{\phi\in
	H^{1/2}(e)}\frac{|\lambda|_{(H^{1/2-\varepsilon}(e))'}|\phi|_{1/2-\varepsilon,e}}{|\phi|_{1/2,e}}\\
\le |\lambda|_{(H^{1/2-\varepsilon}(e))'} h_{e}^\varepsilon
\lesssim h_{e}^\varepsilon \varepsilon^{-1}|\lambda|_{(H^{1/2-\varepsilon}_0(e))'}
\lesssim k^{2\varepsilon}\varepsilon^{-1}|\lambda|_{(H^{1/2}_{00}(e))'},
\end{gather*}
where we have used again \eqref{eq:inversenegative}, and \eqref{penultimadual}.
Therefore,
\[
|\lambda|_{-1/2,e}\lesssim k^{2\varepsilon}\varepsilon^{-1}|\lambda|_{(H^{1/2}_{00}(e))'}.
\]
By taking $\varepsilon=1/(2\log k)$, since
$k^{2\varepsilon}\varepsilon^{-1}$ becomes $2 \exp(1) \log k$, we obtain~\eqref{eq:wish1}.\end{proof}

\section{Stability and error analysis}\label{sect:theory}

In this section, we prove well posedness of the DG formulation in
Problem~\ref{PbGlob} and estimates of the error in the approximation of the
solution to the continuous problem~\eqref{eq:Poisson}.

\subsection{Well posedness}

We prove the well posedness of Problem \ref{PbGlob} by applying 
\cite[Theorem~2.2]{ern.guermond.04}.

%Theorem 1.1, Section II.1 of Brezzi-Fortin. 
%

\newcommand{\hPhi}{\widehat \Phi}

In order to do so, we specify the norms on the discrete spaces:
\begin{align*}
V_h\ &\text{is endowed with the norm}\  \| \cdot \|_{1,\Tess},\\
\Lambda_h\ &\text{is endowed with the norm}\  \tb \cdot \tb_{-1/2,\Squel},
%\\
%\Phi_h\ &\text{is endowed with the norm}\  \| \cdot \|_{1/2,\Squel},
\end{align*}
where the norms $\| \cdot \|_{1,\Tess}$ and $\tb \cdot \tb_{-1/2,\Squel}$
are defined in Section~\ref{sec:globalnorms}.
%and
%\[
%\| \phi \|_{1/2,\Squel}^2 = \sum_{e\subset \Squel}  \| \phi \|^2_{1/2,e}
%\qquad\forall\phi\in \hPhi  =\{
%\phi \in L^2(\Squel): 
%\ \phi|_e \in H^{1/2}(e), \text{for all edge } e \subset \Squel
%\}.
%\]
We also introduce the space $\mathbb{V}_h=V_h\times\Lambda_h$
endowed with the product norm, which is denoted by  $\| \cdot \|_{\mathbb{V}_h}$.

%In order to do so, we introduce the space $\mathbb{V}_h=V_h\times\Lambda_h$
%%$\mathbb{V}_h=  \prod_K \mathbb{P}_k(K) \times \prod_K \Lambda^K$ 
We rewrite Problem \ref{PbGlob} as follows: find $(u,\lambda) \in
\mathbb{V}_h$  such that, for all $(v,\mu) \in
\mathbb{V}_h$, it holds that
\begin{gather}\label{global_discrete}
a(u,\lambda; v,\mu) = F(v,\mu), 
\end{gather}
where 
\begin{multline*}
a(u,\lambda;v,\mu) = \sum_{K\in \Tess} \int_K \nabla u^K \cdot \nabla v^K - \int_{\Squel} \lambda \Lbrack v \Rbrack \cdot \normal + \int_{\Squel} \mu \Lbrack u \Rbrack \cdot \normal\\ + \alpha \sum_K s_K  (D u^K-\gamma_K^* \lambda; tDv^K - \gamma_K^* \mu),
\end{multline*}
and
\begin{equation*}
F(v,\mu) = \int_\Omega f v + \int_{\partial \Omega} g \mu + \alpha \sum_K s_K  (f, tDv^K - \gamma_K^* \mu).
\end{equation*}
%with
%\[
%\begin{split}
%a_K(u^K,\lambda;v^K,\mu) &=  \int_K \nabla u^K \cdot\nabla v^K -
%\int_{\partial K} \lambda (\normalK\cdot \normal) v^K \\ 
%&\qquad\qquad\qquad\qquad+ 
%\int_{\partial K} u^K \mu (\normalK\cdot \normal)  + 
%\alpha s_K  (Du^K - \gamma_K^* \lambda, tDv^K - \gamma_K^* \mu), \\
%F_K(v^K,\mu)&=\int_K f v^K+\alpha s_K  (f, tDv^K - \gamma_K^* \mu).
%\end{split}
%\]

In order to prove the well posedness of Problem \eqref{global_discrete},
we need to prove continuity of the bilinear form
$a(\cdot,\cdot)$  and of the linear functional 
$F(\cdot)$,
as well as an inf-sup condition for $a(\cdot,\cdot)$ in $\mathbb{V}_h \times \mathbb{V}_h$.

\

Remark that, for $u \in H^1(K)$, we have
\begin{equation}\label{LimD}
\| Du \|_{-1,K} = \sup_{0\ne v \in H^1(K)} \frac{\int_K \nabla u \cdot \nabla v}{\| v \|_{1,K}} \leq | u |_{1,K}.
\end{equation}
Then, the following continuity property for the linear functional $F(\cdot)$ is not difficult to prove:  
%There exists a positive constant $C>0$ such that, 
%for all
%$(u,\lambda),(v,\mu)\in \mathbb{V}_h$  it holds that
\begin{gather}
\label{contF}
|F(v,\mu)|\lesssim C(f,g; k)
\|(v,\mu)\|_{\mathbb{V}_h}\qquad\forall (v,\mu)\in \mathbb{V}_h,
\end{gather}
where
\[
 C(f,g; k) =\log(k') 
	\inf_{v \in V: v|_{\partial\Omega} = g} \tb v \tb_{1,\Tess}
  + \| f \|_{0,\Omega} +\alpha  M(k) 
\sqrt{\sum_{K} | f |^2_{-1,K} }.
\]
Indeed, for $v \in V$ with $v|_{\partial\Omega} = g$, letting $\widetilde \mu$ denote the function coinciding with $\mu$ on $\partial\Omega$ and vanishing on $\Squel\setminus\partial\Omega$, we have 
\begin{multline}
\int_{\partial\Omega} g\mu = \sum_{e \subset\partial \Omega} \int_e g \mu = \int_{\Squel}\widetilde \mu \lbrack v \rbrack\cdot \normal \cdot \normal\\ \lesssim  \tb v \tb_{1,\Tess} \tb \widetilde \mu \tb_{-1/2,\Squel} 
\lesssim \log(k')   \tb v \tb_{1,\Tess}  
 \tb  \mu \tb_{-1/2,\Squel}, 
\end{multline}
where we used that, thanks to Proposition \ref{prop:lambdalambda} and to Lemma \ref{lem:log}, we have that \[\tb \widetilde \mu \tb_{-1/2,\Squel} \lesssim \log(k') \tb  \mu \tb_{-1/2,\Squel}.\]
The arbitrariness of $v$ yields 
\[\int_{\partial\Omega} g\mu \lesssim \log(k') \inf_{v \in V:v = g\text{ on }\partial \Omega} \tb v \tb_{1,\Tess} \tb \mu \tb_{-1/2,\Squel}.\]

%Now we have 
%\[
%\sum_{e \subset \partial \Omega} | \bar u^e | = \sum_{e \in \Tess} | \tb \bar u^e - \bar u^K
%\
For the continuity of the bilinear form $a(\cdot,\cdot)$, we start
by observing that %Moreover, we observe that
we can write
\begin{multline}
\int_{\Squel} \lambda \Lbrack v \Rbrack\cdot\normal = \sum_{K \in \Tess} \int_{\partial K}\langle \gamma_K^* \lambda,  v^K \rangle = \sum_{K \in \Tess} \langle \gamma_K^* \lambda , v^K - \bar v^K \rangle +  \sum_{K \in \Tess} \langle \gamma_K^* \lambda ,\bar v^K \rangle
\\ \lesssim \sum_{\in \Tess} | \gamma_K^* \lambda |_{-1,K} | v |_{-1,K} +  \sum_{K \in \Tess}| \bar v^K |\, | \langle \gamma_K^* \lambda , 1 \rangle|.
\end{multline}
If $\lambda$ and $v$ satisfy the following condition
\begin{equation}\label{condlambdav}
\forall K \in \Tess \quad \text{ either }\quad \langle \gamma_K^* \lambda , 1 \rangle = 0 \quad \text{or} \quad \bar v^K = 0,
\end{equation}
then
\[
\int_{\Squel} \lambda \Lbrack v \Rbrack\cdot\normal  \lesssim | v |_{1,\Tess}  \tb \lambda \tb_{-1/2,\Squel}
\]
easily follows.
In the general case, we get the suboptimal bound 
\[
\int_{\Squel} \lambda \Lbrack v \Rbrack\cdot\normal \lesssim  \tb v \tb_{1,\Tess} \tb \lambda \tb_{-1/2,\Squel}.
\]

Therefore, in the general case (without loss of generality we can assume that $\alpha \lesssim 1$), thanks to \eqref{LimD}, using Assumption \ref{sK1}, we have
\begin{multline*}
a(u,\lambda;v,\mu) \lesssim \left(
 \tb u \tb_{1,\Tess}^2 + \tb \lambda \tb_{-1/2,\Squel}^2 
\right)^{1/2}
\left( \tb v \tb_{1,\Tess}^2 
+ \tb \mu \tb_{-1/2,\Squel}^2 
\right)^{1/2}\\ +  M(k) (| u |_{1,\Tess} + \tb \lambda \tb_{-1/2,\Squel})  (| v |_{1,\Tess} + \tb \mu \tb_{-1/2,\Squel}).
\end{multline*}
If $(u,\lambda)$ and $(v,\mu)$ 
 	are such that
 	both $\lambda$ and $v$, and $\mu$ and $u$ satisfy condition \eqref{condlambdav}, then we have
\begin{multline*}
a(u,\lambda;v,\mu) \lesssim | u |_{1,\Tess} | v |_{1,\Tess} + | u |_{1,\Tess} \tb \mu \tb_{-1/2,\Squel} + | v |_{1,\Tess} \tb \lambda \tb_{-1/2,\Squel}\\ +  M(k) (| u |_{1,\Tess} + \tb \lambda \tb_{-1/2,\Squel})  (| v |_{1,\Tess} + \tb \mu \tb_{-1/2,\Squel}),
\end{multline*}
  that yields
\begin{gather}
	\label{contkerb}
	|a(u,\lambda;v,\mu)|\le C\,
        \left(1+M(k)\right)\|(u,\lambda)\|_{\mathbb{V}_h}\|(v,\mu)\|_{\mathbb{V}_h}.
    \end{gather}

We prove the following proposition (inf-sup condition).

\begin{proposition} {\rm (Inf-sup for $a(\cdot,\cdot)$)} We have
\begin{equation}\label{infsupkerb}	\inf_{(u,\lambda) \in \mathbb{V}_h}\sup_{(v,\mu) \in \mathbb{V}_h}
\frac{a(u,\lambda;v,\mu)}{\| (u,\lambda) \|_{\mathbb{V}_h} \| (v,\mu)
	\|_{\mathbb{V}_h}} \gtrsim \left(\!\frac{\rho(k)}{M(k)}\!\right)^2,
\end{equation}
where $M(k)$ and $\rho(k)$ are as in Assumptions~\ref{sK1}
and~\ref{sK2}, respectively.
\end{proposition}

\begin{proof}
 Let $(u,\lambda) \in \mathbb{V}_h$, and let 
\[
v = u -  
\hat v, \quad \text{ with } \quad  \hat v{|_K} =\hat v^K, %(\hat v^K)_K,
\quad \hat v^K =
\int_{\partial K} \lambda (\normalK \cdot \normal) = \langle
\gamma_K^* \lambda, 1 \rangle,
\quad K\in\Tess,
\]
and 
\[
\mu = \lambda + \beta
\hat \mu,\quad \text{ with }\quad 
\quad \hat \mu|_e
  = h_e^{-1}\Lbrack \bar
  u\Rbrack\cdot\normal,\quad e \subset \Squel,
\]
where $\bar u$ denote the piecewise constant function that assumes on each
$K$ the value $\bar u^K$  of the average of $u$ on $K$, and where $\beta$ is a positive constant whose choice will be made later on.  
We observe that 
\begin{equation*}%\label{boundhatv}
\| \hat v \|_{1,*}^2 = 
\sum_{e \subset \Squel} | \sum_{K: e \subset K} \langle \gamma_K^* \lambda , 1 \rangle  |^2 \lesssim \sum_K | \langle \gamma_K^* \lambda, 1 \rangle |^2,
\end{equation*}
where the last bound is obtained by using the fact that any edge belongs to at most two elements, and that the number of edges per element is uniformly bounded by a constant, thanks to Assumption \ref{ass:meshes}.
Moreover, we have that 
\begin{equation}\label{boundhatmu1}
\| \hat \mu \|_{-1/2,e}^2 =  | \Lbrack \bar u \Rbrack |^2,
\end{equation}
which, thanks to Proposition \ref{prop:lambdalambda}, yields
\begin{equation}\label{boundhatmu2}
\tb  \hat \mu \tb^2_{-1/2,\Edges} \lesssim \sum_{e \in \Edges} | \Lbrack \bar u \Rbrack |^2.
\end{equation}
By combining the previous two bounds and applying a triangular inequality, we get
\begin{equation}\label{boundhatV}
\| v,\mu \|_{\mathbb{V}_h} \lesssim \| u,\lambda \|_{\mathbb{V}_h}.
\end{equation}

We can write
\begin{multline*}
a(u,\lambda;v,\mu) = \sum_{K\in\Tess} | u^K |_{1,K}^2 +  \sum_{K\in\Tess} 
{| \langle \gamma_K^* \lambda , 1 \rangle |^2}
 + \beta \sum_{e\subset\Squel} h_e^{-1} \int_{e} \Lbrack u \Rbrack \cdot\Lbrack \bar u \Rbrack  \\
+\alpha \sum_{K\in\Tess} s_K  (Du^K - \gamma_K^* \lambda, t Du^K-\gamma_K^* (\lambda+\beta \hat \mu)).
\end{multline*}

By adding and subtracting $\bar u$, and using a Young inequality, we have
\begin{multline*}
h^{-1}_e	\int_e \Lbrack u \Rbrack \cdot\Lbrack \bar u \Rbrack
=  h^{-1}_e \int_e | \Lbrack \bar u \Rbrack |^2 + h^{-1}_e\int_e
\Lbrack u - \bar u \Rbrack \cdot\Lbrack \bar u \Rbrack
\ge 
\frac{1}{2}| \Lbrack \bar u \Rbrack |^2 
-\frac 1 2 h_e^{-1} \int_e|\Lbrack u - \bar u \Rbrack|^2.
\end{multline*}
We can bound the last term as
\begin{multline*}
\frac 1 2 h_e^{-1} \int_e|\Lbrack u - \bar u \Rbrack|^2
\le h_e^{-1}\left(\int_e| u^{K^+} - \bar u^{K^+}|^2+\int_e|u^{K^-} -
\bar u^{K^-} |^2\right)\\
\le h_e^{-1}\left(\|u^{K^+} - \bar u^{K^+}\|_{0,\partial K^+}^2
+\|u^{K^-} - \bar u^{K^-}\|_{0,\partial K^-}^2\right) 
\lesssim |u^{K^+}|_{1, K^+}^2+|u^{K^-}|_{1, K^-}^2,
\end{multline*}
where in the last inequality we have used the trace inequality \eqref{trace}, the
shape regularity assumption Assumption \ref{ass:meshes}, (ii), and the Poincar\'e inequality \eqref{poincare}.
Therefore, we obtain, 
with some positive constant $c$,
\[
h^{-1}_e	\int_e \Lbrack u \Rbrack \cdot\Lbrack \bar u \Rbrack \ge \frac12|\Lbrack \bar u \Rbrack |^2  - \frac12 c\, (| u^{K^+} |_{1,K^+}^2 + | u^{K^-} |_{1,K^-}^2 ).
\]

% We also observe that
% \begin{equation*}
% \int_{\partial K} \lambda^K \hat v^K \simeq | \partial K |^2 | \bar \lambda^K |^2, \text { with } \bar \lambda^K =  |\partial K|^{-1} \int_{\partial K} \lambda^K.
% \end{equation*}

Then we can write, again for some positive constant $c$,
\begin{multline*}
a(u,\lambda;v,\mu) \geq  \sum_{K\in\Tess} | u^K |_{1,K}^2 +  \sum_{K\in\Tess}
{| \langle \gamma_K^* \lambda , 1 \rangle |^2}  +\frac \beta 2 \sum_{e\subset\Squel}  |
\Lbrack \bar u \Rbrack |^2 \\ - c\beta \sum_{K\in\Tess}  | u^K |_{1,K}^2 +
\alpha \sum_{K\in\Tess} s_K  (Du^K - \gamma_K^* \lambda, t Du^K-\gamma_K^* (\lambda+\beta \hat \mu)) \\  
=\sum_{K\in\Tess} | u^K |_{1,K}^2 
+  \sum_{K\in\Tess}
{| \langle \gamma_K^* \lambda , 1 \rangle |^2}	+ \frac \beta 2 \sum_{e\subset\Squel}  | \Lbrack \bar u \Rbrack |^2 
- c\beta \sum_{K\in\Tess}  | u^K |_{1,K}^2
+\alpha\sum_{K\in\Tess} t \,s_K  (Du^K,Du^K) \\
- \alpha \sum_{K\in\Tess} s_K  (Du^K,\gamma_K^* \lambda) 
- \alpha \beta \sum_{K\in\Tess} s_K   (Du^K,\gamma^*_K \hat \mu)\\
-\alpha\sum_{K\in\Tess} t\, s_K  (\gamma_K^*\lambda,Du^K)
+ \alpha  \sum_{K\in\Tess}  s_K  ( \gamma_K^* \lambda, \gamma_K^* \lambda) 
+\alpha\beta \sum_{K\in\Tess}  s_K  ( \gamma_K^* \lambda, \gamma_K^* \hat \mu ) \\
=(1-c\beta)\sum_{K\in\Tess} | u^K |_{1,K}^2 
+  \sum_{K\in\Tess}
{| \langle \gamma_K^* \lambda , 1 \rangle |^2} 	+ \frac \beta 2 \sum_{e\subset\Squel}  | \Lbrack \bar u \Rbrack |^2 
+T_1+T_2+T_3+T_4+T_5+T_6.
%\geq \\
%	\sum_{K\in\Tess} | u^K |_{1,K}^2 + \frac \beta 2 \sum_{e\subset\Squel}  | \Lbrack \bar u \Rbrack |^2 - c'\beta \sum_{K\in\Tess}  | u^K |_{1,K}^2 + \alpha 	c_1 \sum_{K\in\Tess} \| \gamma_K^* \lambda^K \|_{-1,K}^2\\
%	- \alpha\beta c_2 \| Du^K \|_{-1,K} \| \gamma_K^* \hat \lambda^K \|_{-1,K} - \alpha c_3 \sum_{K\in\Tess} \| Du^K\|_{-1,K} \| \gamma_K^* \lambda^K \|_{-1,K} \\
%	+ c_4 \alpha \beta \| \gamma_K^* \lambda^K \|_{-1,K} \| \gamma_K^* \lambda^K \|_{-1,K}
\end{multline*}
%(here $c$ stands for different possible constants, $c_1$, $c_2$
%etc. are fixed constants). 
We bound separately the terms $T_1$ to $T_6$ on the right-hand side. We immediately observe that $T_5$ yields control on $\gamma_K^* \lambda$. In fact, by Assumption \ref{sK1}, we have
\[
T_5=\alpha\sum_{K\in\Tess}  s_K  ( \gamma_K^* \lambda, \gamma_K^* \lambda ) \geq \alpha\rho(k) \sum_{K\in\Tess} | \gamma_K^* \lambda |_{-1,K}^2.
\]
On the other hand, by Assumption \ref{sK1}, we have
\begin{gather*}
| T_1 |=\left|\alpha\sum_{K\in\Tess} t \,s_K   (Du^K, Du^K) \right| \leq | t |  \alpha M(k)\sum_{K\in\Tess}|u^K|_{1,K}^2,
\end{gather*}
as well as	
\begin{multline*}
|T_2|=\left|\alpha\sum_{K\in\Tess} s_K  (Du^K,\gamma_K^* \lambda) \right| \leq  \alpha M(k) \sum_{K\in\Tess} | u^K |_{1,K} | \gamma_K^* \lambda |_{-1,K}\\ \leq \alpha M(k)\frac{\epsilon_3}{2} \sum_{K\in\Tess} | \gamma_K^* \lambda |^2_{-1,K}  + \alpha M(k)\frac{1}{2 \epsilon_3} \sum_{K\in\Tess} | u^K |_{1,K}^2,
\end{multline*}	
and
\begin{multline*}
|T_4|=\left|\alpha\sum_{K\in\Tess}  t  \,s_K   (\gamma_K^*\lambda, Du^K) \right|\leq
| t | \alpha M(k) \sum_{K\in\Tess} | \gamma_K^* \lambda|_{-1,K} |u^K |_{1,K}\\
\leq | t |\alpha M(k)\frac{\epsilon_1}{2}\sum_{K\in\Tess} | \gamma_K^*
\lambda|_{-1,K}^2
+| t |\alpha M(k)\frac{1}{2 \epsilon_1}\sum_{K\in\Tess} |u^K |_{1,K}^2.
\end{multline*}
Using \eqref{boundhatmu2}, we also have
\begin{multline*}
|T_3|=\left|\alpha\beta\sum_{K\in\Tess} s_K   (Du^K,\gamma^*_K \hat \mu) \right|  \leq  \alpha\beta
M(k)\sum_{K\in\Tess} | u^K |_{1,K} |  \gamma^*_K \hat \mu |_{-1,K} \\\leq
\frac12  \alpha\beta M(k) \sum_{e\subset\Squel} | \Lbrack \bar u \Rbrack |^2 + \frac12  \alpha\beta M(k) \sum_{K\in\Tess} | u^K |_{1,K}^2,
\end{multline*}
and
\begin{multline*}
|T_6|=\left|\alpha\beta\sum_{K\in\Tess}  s_K  ( \gamma_K^* \lambda, \gamma_K^* \hat \mu ) \right| \leq \alpha\beta M(k)\sum_{K\in\Tess} | \gamma_K^* \lambda |_{-1,K} |  \gamma^*_K \hat \mu |_{-1,K}  \\
\leq \alpha\beta M(k) (\sum_{K\in\Tess} | \gamma_K^* \lambda |^2_{-1,K} )^{1/2}  (\sum_{e\subset\Squel} | \Lbrack \bar u \Rbrack |^2)^{1/2} \\ 
\leq \alpha\beta M(k)\frac{\epsilon_2}{2} \sum_{K\in\Tess} | \gamma_K^* \lambda |^2_{-1,K} + \alpha\beta M(k) \frac{1}{2\epsilon_2} \sum_{e\subset\Squel} | \Lbrack \bar u \Rbrack |^2.
\end{multline*}

By combining everything, we obtain
\begin{multline*}	
a(u,\lambda;v,\mu) \geq 	
\left(1-c\beta-\frac{\alpha
	M(k)}{2\epsilon_3}-\frac{\alpha\beta M(k)}{2}-|  t |\alpha
M(k)-\frac{| t |\alpha M(k)}{2\epsilon_1}\right) \sum_{K\in\Tess} | u^K |_{1,K}^2
\\
+ \sum_{K\in\Tess}
{| \langle \gamma_K^* \lambda , 1 \rangle |^2}
+
\frac \beta 2\left(1-\frac{\alpha M(k)}{\epsilon_2}-\alpha M(k)\right) 
\sum_{e\subset\Squel}  | \Lbrack \bar u \Rbrack |^2\\
+\alpha\left(\rho(k)-\frac{\beta M(k)\epsilon_2}{2}-\frac{
	M(k)\epsilon_3}{2}-\frac{|t| M(k)\epsilon_1}{2}\right)
\sum_{K\in\Tess} | \gamma_K^* \lambda |^2_{-1,K}.
\end{multline*}	 

%\COMMENT{Check whether all terms are there.}

We now set $\beta = 1/(2c)$, and we choose $\epsilon_1$,
$\epsilon_2$, and $\epsilon_3$ as $C\rho(k)/M(k)$, with $C$
sufficiently small so that 
$\displaystyle{\left(\rho(k)-\frac{\beta M(k)\epsilon_2}{2}-\frac{
		M(k)\epsilon_3}{2}-\frac{|t|
		M(k)\epsilon_1}{2}\right)\le\frac{\rho(k)}{2}}$.
Recalling that $\rho(k) \le M(k)$, because $\rho(k)$ and $M(k)$ are
coercivity and continuity constants, respectively, 
we choose 
$\alpha=C\rho(k)/(M(k))^2$, with $C$ sufficiently small so that 
all the constants are bounded from below by $\left(\rho(k)/M(k)\right)^2$
(up to a constant).

Observe that neither $\beta$ nor $\alpha$ depend on $h$; $\beta$ is
also independent of $k$, but $\alpha$ depends on $k$ and behaves as
$\rho(k)/(M(k))^2$ for increasing $k$.

With such a choice, for a constant $c_0$ independent of $h$
but dependent on $k$ as $\left(\rho(k)/M(k)\right)^2$, we have
\[
a(u,\lambda;v,\mu) \geq c_0(k) \left(
\sum_{K\in\Tess} | u^K |_{1,K}^2 + \sum_{e\subset\Squel}  | \Lbrack \bar u \Rbrack |^2 + \sum_{K\in\Tess} {| \gamma_K^* \lambda |_{-1,K}^2} +   \sum_{K\in\Tess} | \langle \gamma_K^* \lambda , 1 \rangle |^2
\right).
\]
Therefore, using \eqref{boundhatV}, we conclude that
\begin{gather*}
\sup_{(v,\mu)\in {\mathbb V}_h} \frac{a(u,\lambda;v,\mu)} {\| (v,\mu) \|_{ {\mathbb V}_h}}
\geq \frac{a(u,\lambda; u,\lambda+\beta \hat \lambda)}{\| (u-\hat v,\lambda+\beta \hat \lambda) \|_{ {\mathbb V}_h} }
\gtrsim \left(\!\frac{\rho(k)}{M(k)}\!\right)^2\frac{\| (u, \lambda) \|^2_{ {\mathbb V}_h}}{\| (u, \lambda) \|_{ {\mathbb V}_h}}.
\end{gather*}
\end{proof}

Owing to the continuity properties~\eqref{contF} and~\eqref{contkerb}, and the
  inf-sup condition in Proposition~\ref{infsupkerb}, we apply
  \cite[Theorem 2.2]{ern.guermond.04} and

  conclude with the
following result.

\begin{theorem}\label{th:main} Under Assumptions \ref{ass:meshes}, \ref{sK1} and \ref{sK2},
	Problem \ref{PbGlob} admits a unique solution
        $(\uh,\lah)$. Moreover, the following stability bound for
        $(u_h,\lambda_h)$ is satisfied:
	\[
	\| (\uh,\lah) \|_{\mathbb{V}_h} \lesssim \left(	\frac{M(k)}{\rho(k)}
	\right)^2 \left(\log(k') \inf_{v\in V: v=g\text{ on }\partial\Omega} 	\tb v \tb_{1,\Tess} +
	\| f \|_{0,\Omega} +  \alpha  M(k) 
	\sqrt{\sum_{K} | f |^2_{-1,K} }
	\right).
	\]
\end{theorem}
}

\subsection{Error estimate} We have the following theorem.
\begin{theorem}
	Under Assumptions \ref{ass:meshes}, \ref{sK1} and \ref{sK2}, letting $\w$ denote the solution of \eqref{eq:Poisson}, $\la = \nabla \w \cdot \normal$, and $(\uh,\lah)$ the solution of Problem \ref{PbGlob}, and assuming that $\w \in H^\ell(\Omega)$, $\ell \geq 2$, then the following bound holds:
	\begin{equation}\label{errorestimate}
\| (\w - \uh,\la - \lah) \|_{\mathbb{V}_h} \lesssim (1 + M(k)) \left(
\frac{M(k)}{\rho(k)} 
\right)^2 \frac{h^{s-1}}{k^{\ell-1}} | \w |_{\ell,\Omega},
	\end{equation}
	where $s = \min\{ \ell,k+1\}$.
\end{theorem}

\begin{proof} As we are interested in a $k$-robust estimate, for
    the sake of simplicity we can assume that $k'\not =0$, which is
    always the case except when $k=1$, $k'= k-1$. The latter case, which has little interest in our framework, can be treated with minor modifications to the following arguments.
Let us start by observing that, letting $(\w,\la) \in V \times
\Lambda $ denote the solution to \eqref{eq:PbCont}, for any
$(w,\zeta) \in \mathbb{V}_h$ it holds that
\[a(\w,\la;w,\zeta) = F(w,\zeta).
\]
Then, using \eqref{infsupkerb}, for any $(v,\mu) \in \mathbb{V}_h$ with, for all $K \in \Tess$ and $e \subset \Squel$, $\int_K v^K = \int_K \w$ and $\int_e \mu = \int_e \la$,  we can write:
\begin{multline*}
\left(\frac{\rho(k)}{M(k)}\right)^2\| (\uh - v,\la_h - \mu)
\|_{\mathbb{V}_h} \lesssim \sup_{(w,\zeta) \in \mathbb{V}_h}
\frac{a(\uh-v,\lah - \mu;w,\zeta)}{\| (w,\zeta)
\|_{\mathbb{V}_h}} \\= \sup_{(w,\zeta) \in \mathbb{V}_h} \frac{a(\w-v,\la -
\mu;w,\zeta)}{\| (w,\zeta) \|_{\mathbb{V}_h}} \lesssim (1 + M(k)) \|
(\w-v,\la - \mu)\|_{\mathbb{V}_h}.
\end{multline*}
Using a triangular inequality and the arbitrariness of $\mu$ we obtain
\[
\| (\w-\uh,\la- \lah) \|_{\mathbb{V}_h} \lesssim  (1 + M(k))
\left(\frac{M(k)}{\rho(k)}\right)^2 \inf_{(v,\mu) \in \mathbb{V}_h}  \|
(\w-v,\la - \mu)\|_{\mathbb{V}_h}.
\]

\

In order to bound the right hand side, we recall (see Lemma 23 of
\cite{CangianietalBook}) that, for all $K$ and for all $u \in H^\ell(K)$, with $\ell \geq 1$,
there exist a polynomial $\widetilde \Pi^K_k u \in \mathbb{P}_{k}(K)$
such that, under Assumption \ref{ass:meshes}, we have that
\[
h_K^{-1} \| u - \widetilde \Pi^K_k u \|_{0,K} + 
\| u - \widetilde \Pi^K_k u \|_{1,K} \lesssim \frac{{h_K^{s-1}}}{k^{\ell -
1}} \| u \|_{H^\ell(K)},
\]
where $s = \min\{\ell,k+1\}$,   $\| \cdot \|_{H^\ell(K)}$  denoting the standard, unscaled norm for $H^\ell(K)$:
\[
\| u \|_{H^\ell(K)}^2 = \sum_{j=0}^\ell \sum_{|\alpha| = j} \int_K  \left( \frac{\partial^{|\alpha|} u }{\partial x^\alpha} \right)^2.
\]
 Moreover, for $\lambda \in H^{\ell'}(e)$
 and $\pi_e \lambda$ its $L^2(e)$ projection into $\mathbb{P}_{k'}(e)$,
% the space of degree \rosa{$k'$} polynomials,
 we have
\[
\| \lambda - \pi_e \lambda \|_{0,e} \lesssim{\frac{h_e^{s' }}{(k')^{\ell'}}} |
\lambda |_{\ell,e},
\]
with $s' = \min\{ \ell',k'+1\}$.
Using an Aubin-Nitsche duality argument, we can write
\begin{multline*}
\| \lambda - \pi_e \lambda \|_{-1/2,e} = | \lambda - \pi_e \lambda
|_{-1/2,e} = \sup_{{\phi \in H^{1/2}(e)}:\ {\int_{e} \phi= 0}}
\frac{\int_e (\lambda - \pi_e\lambda) \phi}{| \phi |_{1/2,e}} \\=
\sup_{{\phi \in H^{1/2}(e)}:\ {\int_{e} \phi= 0}} \frac{\int_e (\lambda -
\pi_e\lambda) (\phi - \pi_e \phi)}{| \phi |_{1/2,e}} \lesssim \| \lambda
- \pi_e\lambda \|_{0,e} \sup_{{\phi \in H^{1/2}(e)}:\ {\int_{e} \phi=
0}} \frac{\| \phi - \pi_e \phi \|_{0,e}}{| \phi |_{1/2,e}} \\ \lesssim
\frac{{h_e^{1/2}}}{{(k')^{1/2}}} \| \lambda - \pi_e\lambda \|_{0,e},
\end{multline*}
finally yielding
\[
\| \lambda - \pi_e \lambda \|_{-1/2,e} \lesssim{
\frac{h_e^{s'+1/2}}{(k')^{\ell'+1/2}}} | \lambda |_{\ell',e}.
\]

\

Let now the solution of \eqref{eq:Poisson} satisfy $u \in
H^\ell(\Omega)$ with $\ell \geq 2$. As $\nabla u \in H^{\ell -
1}(\Omega) \subseteq H^1(\Omega)$ we have that for $e$ edge of $K$,
$\nabla u \cdot \normal \in H^{\ell - 3/2}(e) \subseteq
H^{1/2}(e)$ and $\| \nabla u \cdot \normal \|_{\ell - 3/2,e}
\lesssim \| \nabla u \|_{\ell - 1,K}$. Letting then $v^K = \widetilde \Pi_k^K(u)$ and
$\mu|_e = \pi_e (\nabla u \cdot \normal)$,
    and observing that for $k>1$ we have $k' \gtrsim k$, thanks to \eqref{stimanormath1} we have 
\begin{gather*}
\| (u - v,\lambda - \mu) \|^2_{\mathbb{V}_h} \lesssim \frac{h_K^{2(s-1)}}{k^{2(\ell - 1)}} 
\sum_{K\in\Tess}  \left(   \| u \|^2_{\ell,K} + 
\sum_{e \subset \partial K}   \| {\nabla u \cdot \normal} \|^2_{\ell -
3/2,e}\right)
 \lesssim \sum_{K\in\Tess}\frac{h_K^{2(s-1)}}{k^{2(\ell - 1)}}    \| u \|^2_{\ell,K},
%\lesssim \frac{h^{s-1}}{k^{\ell - 1}} \| u
%\|_{\ell,\Omega},
\end{gather*}
which concludes the proof.
\end{proof}

\newcommand{\hlambdaK}{\widehat{\lambda}^K}

\newcommand{\hmu}{\widehat{\mu}}
\newcommand{\hmuK}{\widehat{\mu}^K}
\newcommand{\hlambda}{\widehat \lambda}
\newcommand{\LambdaK}{\Lambda_K}
\newcommand{\hLambda} {\widehat\Lambda_h}
\newcommand{\aK}{\widehat a^K}

\section{Stabilization forms}\label{sect:stab}

%\COMMENT{Anche qui \`e rimasto tutto invariato, tranne che si introduce $\Lambda_K$ che non era stato definito e si sostituisce $k$ con $k'$.}

In order for the proposed method to be practically feasible, we need
to construct computable bilinear forms $s_K(\cdot,\cdot)$ satisfying Assumptions~\ref{sK1}
and~\ref{sK2}. We follow
%This can be done by using
the approach of \cite{SilviaPaperAltro}. Let %Following such paper, and letting
  \begin{equation}\label{defLambdaK}
	\Lambda_K = \Lambda_h |_{\partial K} = \{ \lambda \in L^2(\partial K): \ \lambda|_e \in \mathbb{P}_{k'}(e), \ e \subset \partial K\}\end{equation} 
(we recall that $k'\in \{k,k-1\}$), and introduce an auxiliary
space $W^K \subseteq H^1(K)$ with $W^K \cap \mathbb{P}_0(K) = \{0\}$,
and with $\dim(W^K) = n_K = \dim(\Lambda_K)-1$, satisfying, for some
positive constant $\rho(k')$, an inf-sup condition of the form
\begin{equation}\label{mangiacoda}
\inf_{\lambda \in \Lambda_K: \ \int_{\partial K}\lambda=0 } 
\blu{\sup_{w \in W^K} }
%\sup_{w \in W^K: \int_{\partial K} = 0} 
\frac{\int_K \lambda w}{| \gamma_K^ * \lambda |_{-1,K} | w |_{1,K}} \gtrsim \rho(k).
\end{equation}
The choice of the subspace $W^K$ that characterize our method is
  specified below.

\renewcommand{\s}{\sigma}
\renewcommand{\S}{\Sigma}

Let $\phi_i$, $i=1,\cdots,\blu{n_K}$, denote a basis for $W^K$.
Consider 
the operator $\s:H^1(K) \times H^1(K) \to \mathbb{R}$ given by
\[\s(w,v) = \int_K \nabla w \cdot \nabla v.
\]
We observe that we have
\[
| \s(w,v) | \leq | w |_{1,K} | v|_{1,K},\qquad \s(w,w) = | w |_{1,K}^2.
\]
We also observe that, as $W^K \cap \mathbb{P}_0(K) = \{0\}$, the seminorm $| \cdot |_{1,K}$ is a norm on $W^K$. 
We let $\S$ denote the stiffness matrix associated with the restriction
of $\s(\cdot,\cdot)$ to $W^K$, i.e. 
\[\S_{ij}=\s(\phi_j, \phi_i), \qquad i,j=1,\cdots,\blu{n_K}.\]

\

We can now introduce the following bilinear form $s_K: (H^{1}(K))' \times (H^{1}(K))'\to \mathbb{R}$ defined as
\[
s_K(F,G) = \vec F^T \S^{-1} \vec G, \quad \text{ with } \vec F = (\langle F,\phi_i\rangle)_{i=1}^{n_K}, \  \vec G = (\langle G,\phi_i\rangle)_{i=1}^{n_K}.
\]

It is not difficult to prove that the bilinear form $s_K(\cdot,\cdot)$ satisfies
Assumption \ref{sK1} with $M(k)=1$.
Moreover, it is possible to prove 
(see \cite{
SilviaPaperAltro}) 
that, provided that (\ref{mangiacoda})
holds, $s_K(\cdot,\cdot)$ satisfies also Assumption~\ref{sK2} (actually, (\ref{mangiacoda}) is a necessary and sufficient condition for Assumption~\ref{sK2} to hold).

Observe that, for $u,v \in( H^{1}(K))'$ and $\lambda,\mu \in H^{-1/2}(\partial K)$, we have
\begin{equation}\label{eq:sK-howto}
s_K(\D u - \gamma_K^* \lambda,t \D v - \gamma_K^* \mu) = \vec \eta^T \S^{-1} \vec \zeta, \qquad s_K (f,t \D v - \gamma_K^* \mu ) = \vec f^T \S^{-1} \vec \zeta
\end{equation}
with
\begin{equation}\label{eq:sK-howto-vectors}
\eta_i = \int_K \nabla u \cdot \nabla \phi_i - \int_{\partial K} \lambda \phi_i, \qquad \zeta_i = t \int_K \nabla v \cdot \nabla \phi_i - \int_{\partial K} \mu \phi_i, \qquad f_i = \int_K f \phi_i.
\end{equation}

\

In order to complete the definition of our method, we only need to choose the subspace $W^K$ of $H^1(K)$.
In order to do that, we subdivide the polygonal element $K$ into $N_K$
triangles $T_i$, $1\le i\le N_K$, each having one edge, denoted by $e_i$, coinciding with one edge of $K$, and
the opposite vertex coinciding with $\mathbf{x}_K$, the center of the
ball in Assumption~\ref{ass:meshes}, (i). Due to Assumption~\ref{ass:meshes}, all the triangles $T_i$ are shape
regular.
Let us consider a reference triangle $\widehat{T}$ and denote by $F_{i}$
the affine maps from $\widehat{T}$ to $T_i$, defined in such a way
that the edge $\widehat{e}$ is mapped onto $e_i$.
We construct a finite dimensional space $\widehat{W}\subseteq H^1(\widehat{T})$ as follows.

Let $\widehat{V}_\delta\subseteq H^1(\widehat{T})$ be a family of finite
dimensional approximation spaces, whose elements vanish on
$\partial\widehat{T}\setminus \widehat{e}$, each \blu{constructed on a quasi uniform mesh of $\widehat{T}$} of mesh size $\delta$. \blu{The approximation
  assumptions $\widehat{V}_\delta$ needs to satisfy are stated in
  Lemma~\ref{lem:infsuphat}
and Theorem~\ref{th:infsupK} below; a specific choice will be given in Section \ref{sec:results}}.

We define the operator $\mathcal{G}:\mathbb{P}_{k'}(\widehat{e})\to
\widehat{V}_\delta$ that maps $\widehat{\lambda}\in
\mathbb{P}_{k'}(\widehat{e})$ to the (unique) function $\widehat{\phi}\in
\widehat{V}_\delta$
that satisfies
\begin{equation}\label{eq:defphihat}
\int_{\widehat{T}}\nabla\widehat{\phi}\cdot \nabla\widehat{v}= 
\int_{\widehat{e}}\widehat{\lambda}\widehat{v}\qquad
\text{ for all } \widehat{v}\in \widehat{V}_\delta.
\end{equation}
Notice that $\mathcal{G}(\hat\lambda)$ is a discretized harmonic
  lifting in $\hat T$ of the Neumann datum $\hat\lambda$ on $\hat e$; see the proof of
Lemma~\ref{lem:infsuphat} below.
We then define $\widehat{W}$ as
\[
\widehat{W}=\mathcal{G}(\mathbb{P}_{k'}(\widehat{e})).
\]
We set
\begin{equation}\label{eq:defWi}
W_{i}=\{\widehat{w}\circ F_i^{-1} \text{ with } \widehat{w}\in
\widehat{W}\},
\end{equation}
and
\begin{equation}\label{eq:defWK}
W^K=\{w\in L^2(K):\; w|_{T_i}\in W_{i},\ 1\le i\le N_K\}.
\end{equation}
Notice that, as the functions in $W_{i}$ have zero Dirichlet traces
along the edges of each $T_i$  interior to $K$, the functions of $W^K$ are continuous; therefore
$W^K\subset H^1(K)$.  We also remark that, in order to construct
  $W^K$ for any $K\in\Tess$, one needs to solve~\eqref{eq:defphihat}
  for each function of a basis of $\mathbb{P}_{k'}(\widehat{e})$ on
  the reference element. This can be done  offline once and for all; for more details, see Section \ref{sec:results}.

We prove the following inf-sup condition on the reference triangle
$\widehat{T}$.

\begin{lemma}\label{lem:infsuphat}
Assume that the space $\widehat{V}_\delta$ is such that, for all
$v\in H^2(\widehat{T})$ with $v=0$ on $\partial\widehat{T}\setminus
\widehat{e}$,
it holds
\begin{equation}\label{approxVdelta}
\inf_{v_\delta\in \widehat{V}_\delta}\|v-v_\delta\|_{1,\widehat{T}}\lesssim \delta
|v|_{2,\widehat{T}}.
\end{equation}
Then there exists a constant $c_0>0$ independent of $k'$ such that, provided
that $\delta<c_0(k')^{-2}$, we have
\[
\inf_{\lambda\in \mathbb{P}_{k'}(\widehat{e})}\sup_{\widehat{\phi}\in \widehat{W}}
\frac{\int_{\widehat{e}}\lambda
  \widehat{\phi}}{\|\lambda\|_{(H^{1/2}_{00}(\widehat e))'} %\|\lambda\|_{(H^{1/2}_{00}(\widehat{e}))'}
|\widehat{\phi}|_{1,\widehat{T}}}\gtrsim 1.
\]
\end{lemma}

\begin{proof}
Fix $\lambda\in \mathbb{P}_{k'}(\widehat{e})$. We let $u_\lambda$ denote the solution to
\begin{equation}\label{auxpbstrong}
-\Delta u_\lambda = 0 \ \text{in }\widehat{T}, \qquad u_\lambda = 0\ \text{on }\partial\widehat{T}\setminus \widehat e, \qquad \nabla u_\lambda \cdot \mathbf{n} = \lambda \ \text{on }\widehat e.
\end{equation}
Writing \eqref{auxpbstrong} in variational form, we easily see that $u_\lambda$ satisfies
\begin{equation}\label{eq:auxpbl}
\int_{\widehat{T}}\nabla u_\lambda\cdot\nabla
v=\int_{\widehat{e}}\lambda v \qquad \text{ for all } v\in
H^1(\widehat{T}) \text{ with } v=0 \text{ on } \partial\widehat{T}\setminus \widehat{e}.
\end{equation}

We can then write
\[
\int_{\widehat{e}}u_\lambda \lambda
=\int_{\widehat{T}}|\nabla u_\lambda|^2=|u_\lambda|_{1,
  \widehat{T}}^2\gtrsim
|u_\lambda|_{1, \widehat{T}} \|\lambda\|_{(H^{1/2}_{00}(\widehat e))'},
\]
where, using~\eqref{auxpbstrong}, the last bound follows from
\begin{multline}\label{eq:boundLambda}
\|\lambda\|_{(H^{1/2}_{00}(\widehat e))'}
=\sup_{\phi\in
	H^{1/2}_{00}(\widehat{e})}
\frac{\int_{\widehat{e}}\lambda\phi}{\|\phi\|_{H^{1/2}_{00}(\widehat{e})}}
\simeq \sup_{\phi\in{H^{1/2}(\partial\widehat T),\
		\phi_{|_{\partial\widehat T\setminus \widehat{e}}}=0}}
\frac{\int_{{\widehat e}}\lambda\phi}{|\phi|_{H^{1/2}(\partial\widehat{T})}} \\ = 
\sup_{\phi\in{H^{1/2}(\partial\widehat T),\
		\phi_{|_{\partial\widehat T\setminus \widehat{e}}}=0}}
\frac{\int_{\widehat T} \nabla u_\lambda\cdot \nabla \widetilde\phi}{|\phi|_{H^{1/2}(\partial\widehat{T})}}  \lesssim | u_\lambda |_{1,\widehat T},
\end{multline}
$\widetilde \phi \in H^1(\widehat  T)$ denoting the harmonic lifting of $\phi$.

Let now $\widehat{\phi}=\mathcal{G}(\lambda)\in \widehat{W}$.
It is easily seen that $\widehat \phi$ is the Galerkin projection of $u_\lambda$
onto $\widehat V_\delta$. Then it holds
\begin{equation}\label{eq:boundLambda1}
|\widehat{\phi}|_{1,\widehat{T}}\le |u_\lambda|_{1,\widehat{T}}.
\end{equation}
\blu{By using \eqref{eq:inversenegativemezzo},   }
we have
\begin{multline}\label{eq:boundLambda0}
\int_{\widehat{e}}\lambda \widehat{\phi}=\int_{\widehat{e}}\lambda(\widehat{\phi}-u_\lambda)+
\int_{\widehat{e}}\lambda u_\lambda
\ge
-\|\lambda\|_{0,\widehat{e}} \|u_\lambda-\widehat{\phi} \|_{0,\partial\widehat{T}}
+|u_\lambda|_{1, \widehat{T}}
\|\lambda\|_{(H^{1/2}_{00}(\widehat e))'}\\
\blu{\ge} -k' \|\lambda\|_{(H^{1/2}_{00}(\widehat e))'}\|u_\lambda-\widehat{\phi}
\|_{0,\partial\widehat{T}}
+|u_\lambda|_{1, \widehat{T}}
\|\lambda\|_{(H^{1/2}_{00}(\widehat e))'}\\
\gtrsim -k' \|\lambda\|_{(H^{1/2}_{00}(\widehat e))'}
\|u_\lambda-\widehat{\phi}\|_{0,\widehat{T}}^{1/2}
\|u_\lambda-\widehat{\phi}\|_{1,\widehat{T}}^{1/2}+|u_\lambda|_{1, \widehat{T}}
\|\lambda\|_{(H^{1/2}_{00}(\widehat e))'}\\
\gtrsim -\delta^{1/2}k'
\|\lambda\|_{(H^{1/2}_{00}(\widehat e))'}
|u_\lambda|_{1,\widehat{T}}+|u_\lambda|_{1, \widehat{T}}
\|\lambda\|_{(H^{1/2}_{00}(\widehat e))'},
\end{multline}
where we have used the trace bound \eqref{trace}, the Aubin-Nitsche
duality argument, which, thanks to \eqref{approxVdelta}, allows us to estimate $\|u_\lambda - \widehat \phi\|_{0,\widehat T}$ with $\delta | u_\lambda - \widehat \phi |_{1,\widehat{T}}$, and \eqref{eq:boundLambda1}. 
Then we have
\begin{equation*}
\int_{\widehat{e}}\lambda \widehat{\phi} \ge | \widehat \phi |_{1,\widehat{T}} \| \lambda \|_{(H^{1/2}_{00}(\widehat e))'} (c_1 - c_2 \delta^{1/2}k' ).
\end{equation*}

Therefore, by choosing $\delta$ in such a way that
  $(c_1 - c_2 \delta^{1/2}k' ))>0$ (i.e. $\delta < c_1 /(c_2(k')^{2})$), and 
inserting~\eqref{eq:boundLambda1} into~\eqref{eq:boundLambda0}, after
dividing by $|\widehat\phi|_{1, \widehat{T}}$, %and taking $\delta$ such that
%$\delta^{1/2}k$ is small enough 
we get the thesis.
\end{proof}

Due to a scaling argument, thanks to the shape regularity of the
triangles $T_i$, the inf-sup condition of Lemma~\ref{lem:infsuphat} implies 
the following inf-sup condition on $T_i$:
\begin{equation}\label{eq:infsupTi}
\inf_{\lambda\in \mathbb{P}_k(e)}\sup_{{\phi}\in {W_i}}
\frac{\int_{{e}}\lambda
  \phi}{\|\lambda\|_{(H^{1/2}_{00}({e}))'}|\phi|_{1,{T_i}}}\gtrsim 1,
\end{equation}
provided
that the parameter $\delta$ entering the definition of $\widehat
V_\delta$ in the construction of $\widehat W$ satisfies $\delta<c_0(k')^{-2}$, with $c_0>0$
given in Lemma~\ref{lem:infsuphat}.

By mapping and assembling on all subtriangles $T_i$ of $K$, we obtain the
following inf-sup condition on $K$.

\begin{theorem}\label{th:infsupK}
 Let the
space $W^K$ be defined as in \eqref{eq:defWK} with $\delta<c_0(k')^{-2}$, where the constant $c_0>0$
is given in Lemma~\ref{lem:infsuphat}.
Then, under Assumption~\ref{ass:meshes}, we have, for all $K \in \Tess$,
\[
\inf_{\lambda\in \Lambda^K}\sup_{{\phi}\in W^K}
\frac{\int_{\partial K}\lambda
  \phi}{|\lambda|_{-1/2,\partial K}|\phi|_{1,K}}\gtrsim (\log k')^{-1}.
\]
\end{theorem}

\begin{proof}
Let $\{e_i\}_{i=1}^{N_K}$ denote the set of edges of $K$.
Recall that $\Lambda^K=\prod_{1\le i\le N_K}\mathbb{P}_{k'}(e_i)$ and 
$W^K \sim\Pi_{1\le i\le N_K} W_{i}$.

By a standard
argument as in~\cite{SilviaAngela}, from the local inf-sup conditions~\eqref{eq:infsupTi}, we have
\[
\inf_{\lambda\in\Lambda^K}\sup_{\phi\in W^K}
\frac{\sum_{i=1}^{N_K}\int_{e_i}\lambda\phi}{\left(\sum_{i=1}^{N_K}
    |\lambda|_{(H^{1/2}_{00}(e_i))'}^2\right)^{1/2}
  \left(\sum_{i=1}^{N_K} |\phi|_{1,T_i}^2\right)^{1/2}}
\gtrsim 1.
\]
As $\sum_{i=1}^{N_K}\int_{e_i}\lambda\phi=\int_{\partial K}
\lambda\phi$
and $\left(\sum_{i=1}^{N_K}
  |\phi|_{1,T_i}^2\right)^{1/2}=|\phi|_{1,k}$, we only need to prove
the bound
\begin{equation}\label{eq:wish}
|\lambda|_{-1/2,\partial K}
=
\sup_{\phi\in H^{1/2}(\partial K),\ \int_{\partial
    K}\phi=0}\frac{\int_{\partial K}\lambda\phi}{|\phi|_{1/2,\partial K}}
\lesssim 
\log k'
\left(\sum_{i=1}^{N_K}
    |\lambda|_{(H^{1/2}_{00}(e_i))'}^2\right)^{1/2}.
\end{equation}

From the Cauchy-Schwarz inequality, we have
\[
\int_{\partial K}\lambda\phi=\sum_{i=1}^{N_K}\int_{e_i}\lambda\phi
\le \sum_{i=1}^{N_K}\|\lambda\|_{-1/2,e_i}\|\phi\|_{1/2,e_i}
\le \left(\sum_{i=1}^{N_K}\|\lambda\|_{-1/2,e_i}^2\right)^{1/2}
\left(\sum_{i=1}^{N_K}\|\phi\|_{1/2,e_i}^2\right)^{1/2}.
\]
Now, on the one hand, denoting by $\bar\phi^{e_i}$ the average of $\phi$ on
$e_i$, i.e. $\bar\phi^{e_i}=|e_i|^{-1}\int_{e_i}\phi$,
we have
\[
\sum_{i=1}^{N_K}\|\phi\|_{1/2,e_i}^2=
\sum_{i=1}^{N_K}\left( |\bar\phi^{e_i}|^2+ |\phi|_{1/2,e_i}^2\right)
\le
\sum_{i=1}^{N_K}|\bar\phi^{e_i}|^2 + |\phi|_{1/2,\partial K}^2.
\]
For the first term on the righ-hand side, as $\phi$ has zero mean
value on $\partial K$, we obtain
\[
\sum_{i=1}^{N_K}|\bar\phi^{e_i}|^2= \sum_{i=1}^{N_K}
|e_i|^{-2}\left|\int_{e_i} \phi\right|^2\le 
\sum_{i=1}^{N_K}|e_i|^{-1}\int_{e_i}|\phi|^2\lesssim
h_K^{-1} \int_{\partial K}|\phi|^2\lesssim |\phi|_{1/2, \partial K}^2,
\]
where we also have used $|e_i|\simeq h_K$ (due to shape regularity),
and the Poincar\'e inequality \eqref{poincarev2}.
Therefore, 
\[
\frac{\int_{\partial K}\lambda\phi}{|\phi|_{1/2, \partial K}}\lesssim 
\left(\sum_{i=1}^{N_K}\|\lambda\|_{-1/2,e_i}^2\right)^{1/2}.
\]
In order to prove~\eqref{eq:wish}, we only need to apply inequality \eqref{eq:wish1}

% |\lambda|_{-1/2,\partial K} = \sup_{\phi\in H^{1/2}(\partial
%   K)}\frac{\int_{\partial K}\lambda\phi}{|\phi|_{1/2,\partial K}}
% \]
\end{proof}
As $k' \leq k$, Theorem \ref{th:infsupK} yields \eqref{mangiacoda}
with $\rho(k) = (\log k)^{-1}$.

\begin{remark}\label{rem:latini}
 Assumption \ref{ass:meshes} (ii) is needed in the proof of Theorem
 \ref{th:infsupK}, as, under such an assumption, we manage to bound
 the $H^{-1/2}(\partial K)$ semi norm of $\lambda$ with the sum of its
 $(H^{1/2}_{00}(e_i))'$ norms over all edges $e_i$ of $\partial K$. 
% This bound \ip{is} not generally \ip{valid} if $K$ has very small
% edges, so that, in itself, the stabilization that we propose here does
% not allow robustness with respect to decreasing edge length.
This bound is not generally valid if $K$ has very small
edges. Therefore,
% so that, in itself,
the stabilization that we propose here is not proven to be robust
%does not allow robustness
with respect to decreasing edge length. 

However we are confident that suitably combining the present approach
with the approach used in \cite{SilviaDaniele} will allow us
to obtain a method that is simultaneously robust with respect to decreasing edge length and increasing polynomial degree.
\end{remark}

\section{Hybridization}\label{sec:hybrid}
	As in \cite{SilviaDaniele}, in order to efficiently implement the method, we perform an hybridization procedure by introducing an auxiliary unknown $\phi$ approximating the trace on $\Squel$ of the solution and by using independent unknowns $\hlambda^K\in \Lambda_K$ ($\Lambda_K$ defined by \eqref{defLambdaK}) to approximate  $\nabla  \w \cdot \normalK$. To this aim, we introduce the following discrete spaces
	\begin{gather*} \hLambda = \prod_{K\in\Tess} \LambdaK,\qquad
	\Phi_h =\{ \phi \in L^2(\Squel): \ \phi|_e \in \mathbb{P}_{k'}(e), \ e \subset \Squel\} .
	\end{gather*}
	Letting $b: \hLambda \times \Phi_h\to \mathbb{R}$ be defined by
	\[
	b(\hlambda,\psi) = \sum_K \int_{\partial K } \hlambda^K \psi,
	\]
	it is easy to check that $\Lambda_h$ is isomorphic to 
	\[
	\ker b = \{\hlambda \in \hLambda: \ b(\lambda,\psi) =0 \ \forall \psi\in \Phi_h\} \subset \hLambda.
	\]
	More precisely, $\hlambda \in \ker b$ if and only if $\hlambda^K = \lambda (\normal\cdot\normalK)$ for some $\lambda \in \Lambda_h$. Introducing the bilinear forms $\aK : \mathbb{P}_k(K) \times \LambdaK\to \mathbb{R}$ given by
	\begin{multline*}\aK(u^K,\hlambda^K;v^K \hmuK  )  =  \int_K \nabla u^K \cdot \nabla v^K - 
	\int_{\partial K} \hlambda^K v^K + \int_{\partial K} \hmuK u^K
	\\ + \alpha  s_K (D u^K-\gamma_K^* \hlambda^K; tDv^K - \gamma_K^* \hmuK),
	\end{multline*}
	and letting
	\[
	\widehat a(u,\hlambda;v,\hmu) = \sum_K \aK(u^K,\hlambda^K;v^K,\hmuK), \qquad F(v,\hmu) = \int_\Omega f v + \int_{\partial\Omega} g\mu + \alpha \sum_{K\in \Tess} s_K (f; tD v^K - \gamma_K^* \mu),
	\]
	we can then consider the following hybridized problem.
	\begin{problem}\label{PbHybrid} Find $\uh = (\uh^K)_{K\in\Tess} \in V_h$, $\hlah = (\hlah^K) \in \hLambda$, $\phi \in \Phi_h$ with $\phi|_{\partial\Omega} = g$ such that, for all $v = (v^K)_{K\in\Tess} \in V_h$, $\hmu = (\hmu^K) \in \hLambda$, $\psi \in \Phi_h$ with $\psi|_{\Omega} = 0$,
		\begin{gather}
		\label{localpb}	\widehat a(\uh,\hlah; v,\hmu) - b(\hmu,\phi) = F(v,\hmu)\\
\label{eq:PbGlob-glue}		b(\hlah,\psi) = 0.
		\end{gather}
		
	\end{problem}
	
	The well posedness of Problem \ref{PbHybrid} and its equivalence to Problem \ref{PbGlob} are proven in \cite{SilviaDaniele}. Observe that \eqref{localpb} reduces to independent Dirichlet problems in each $K$, with boundary condition $u^K = \phi$ on $\partial K$, and with non standard stabilization given by the bilinear form $s_K $. The local unknown can then be eliminated by static condensation, reducing the solution to a problem on the unknown $\phi$.

\section{Numerical Results}\label{sec:results}

\begin{figure}
\centering
\subfloat{\includegraphics[width=0.3\textwidth]{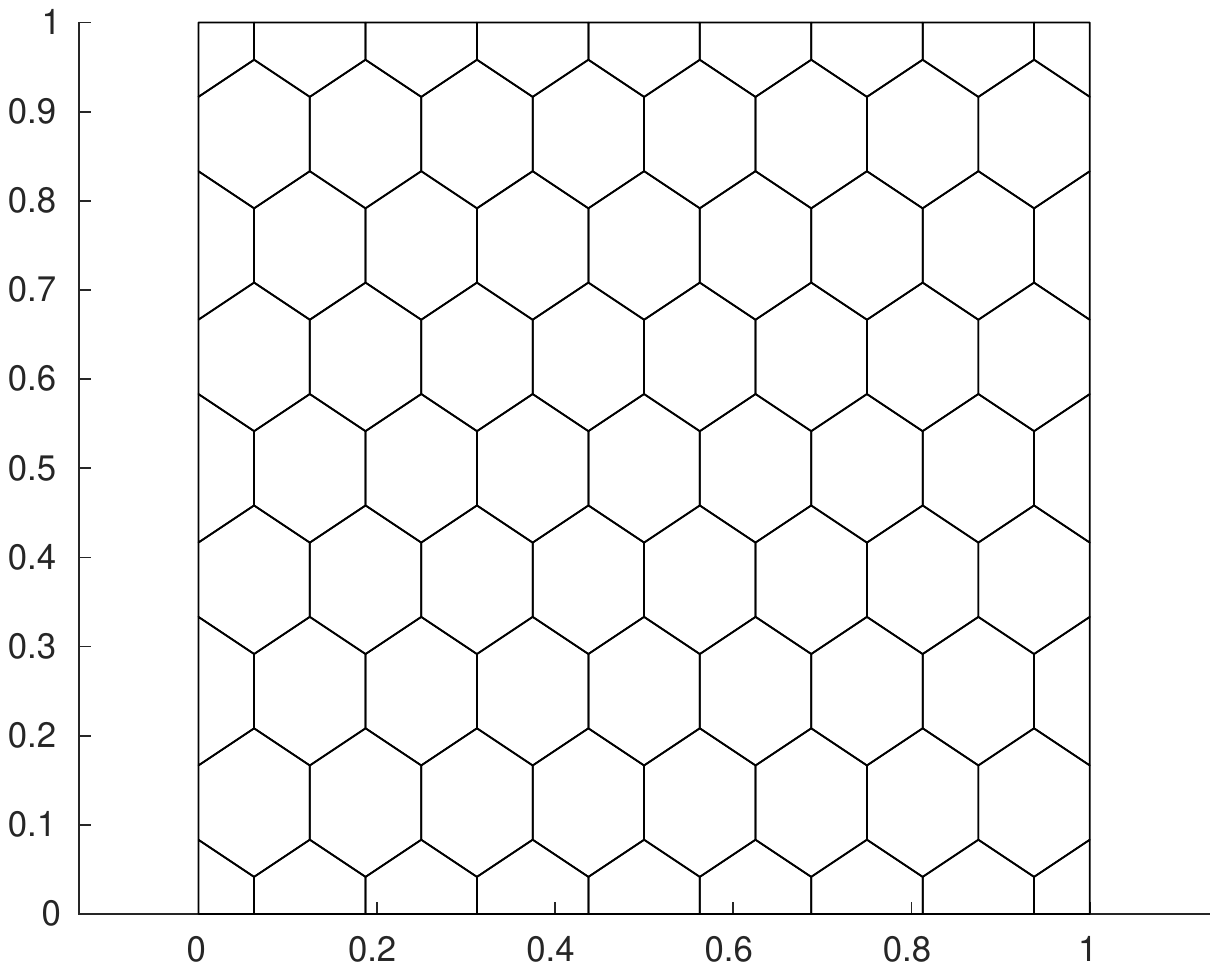}\label{fig:rhexa}}
\subfloat{\includegraphics[width=0.3\textwidth]{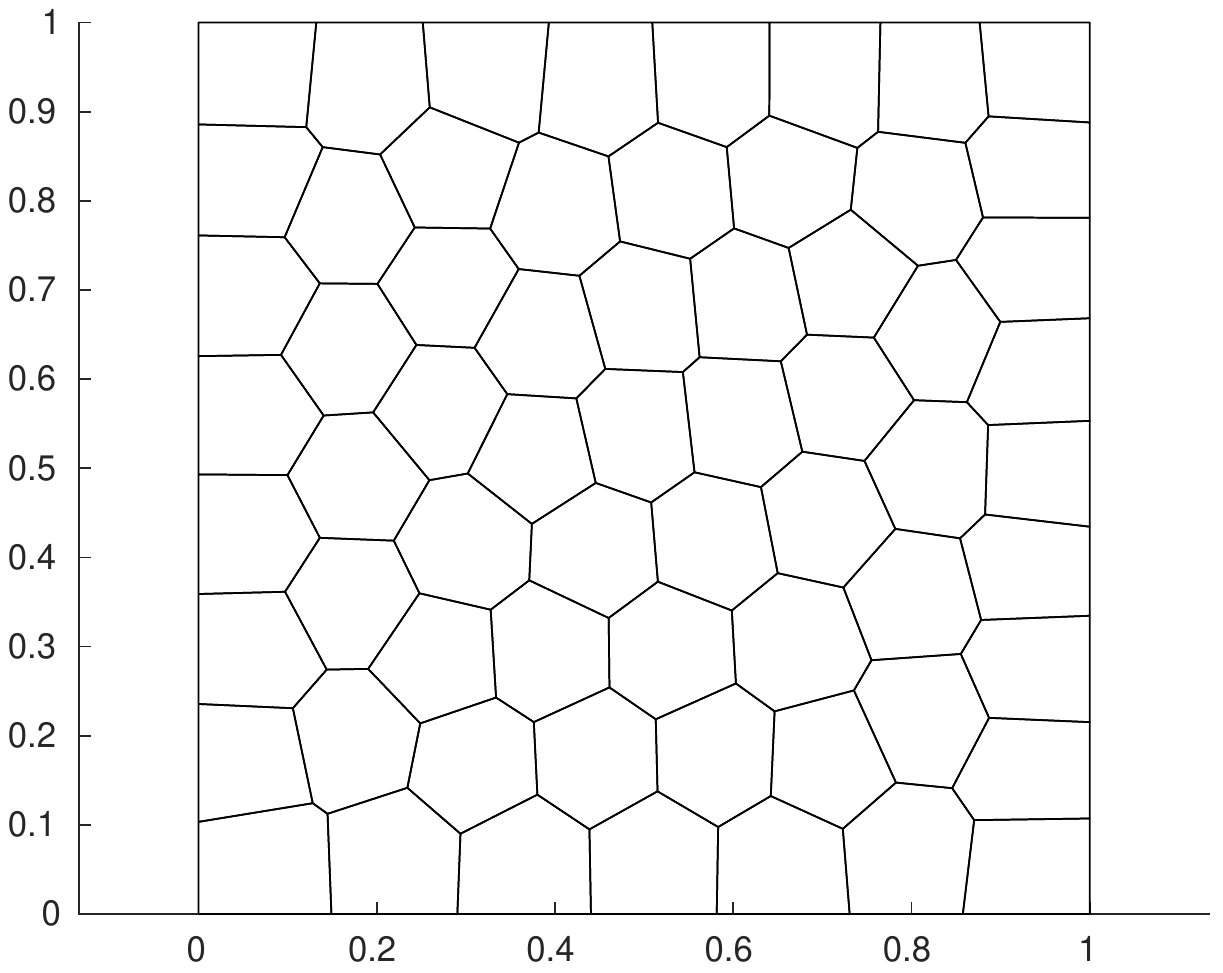}\label{fig:cvt}}
\subfloat{\includegraphics[width=0.3\textwidth]{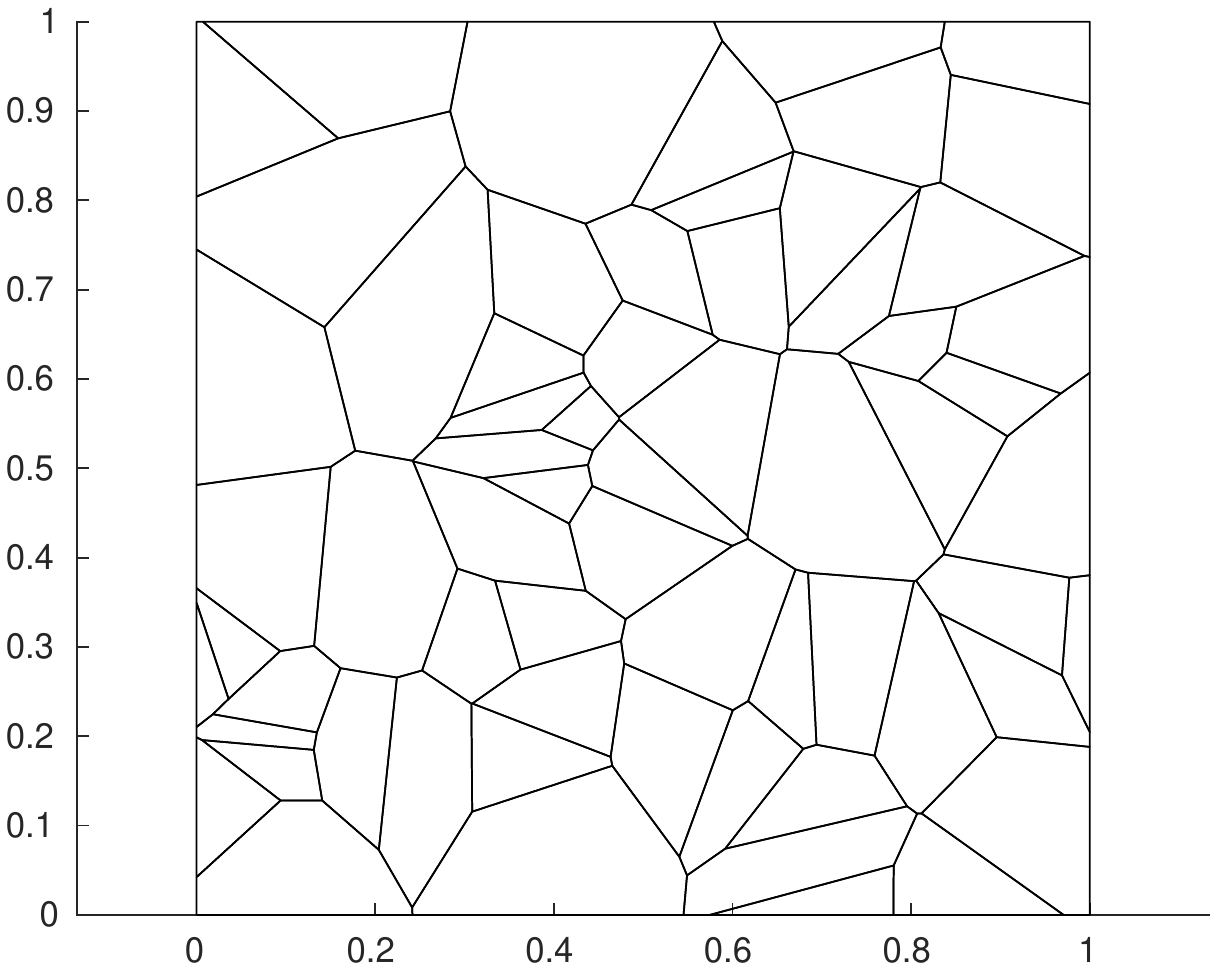}\label{fig:voro}}
\caption{Meshes used in experiment ii). From left to right: meshes made of regular hexagons, Central Voronoi Tessellation, and random Voronoi cells.}
\label{fig:meshes-k-robustness}
\end{figure}

\begin{table}
\centering
\caption{Meshes of regular hexagons used in experiment i). For all these meshes, $\gamma_0 \approx \num[scientific-notation=fixed, fixed-exponent=0]{5.33}, \gamma_1 \approx \num[scientific-notation=fixed, fixed-exponent=0]{3.16}$.}
\label{tab:mesh-r-hexa}
\begin{tabular}{
c
S[table-format=5.0]
S[table-format=6.0]
S[table-format=1.{\roundPrecision}e-1]
S[table-format=1.{\roundPrecision}e-1]
}
\toprule
Mesh & {$N_p$} & {$N_e$} & {$h$} & {$h_\textup{min}$}\\
\midrule
r-hexa$_{1}$ & 4193 & 12580 & 2.083333e-02 & 5.208333e-03\\
r-hexa$_{2}$ & 10151 & 30454 & 1.333333e-02 & 3.333333e-03\\
r-hexa$_{3}$ & 22726 & 68179 & 8.888889e-03 & 2.222222e-03\\
r-hexa$_{4}$ & 40301 & 120904 & 6.666667e-03 & 1.666667e-03\\
r-hexa$_{5}$ & 62876 & 188629 & 5.333333e-03 & 1.333333e-03\\
r-hexa$_{6}$ & 90451 & 271354 & 4.444444e-03 & 1.111111e-03\\
\bottomrule
\end{tabular}
\end{table}

\

\begin{table}
\centering
\caption{Central Voronoi Tessellations used in experiment i).}
\label{tab:mesh-cvt}
\begin{tabular}{
c
S[table-format=5.0]
S[table-format=6.0]
S[table-format=1.{\roundPrecision}e-1]
S[table-format=1.{\roundPrecision}e-1]
S[table-format=1.2]
S[table-format=1.{\roundPrecision}e+1]
}
\toprule
Mesh & {$N_p$} & {$N_e$} & {$h$} & {$h_\textup{min}$} & {$\gamma_0$} & {$\gamma_1$}\\
\midrule
cvt$_{1}$ & 2048 & 6124 & 3.400609e-02 & 1.506563e-03 & 3.181507 & 1.878109e+01\\
cvt$_{2}$ & 4096 & 12250 & 2.698981e-02 & 1.008020e-03 & 3.504224 & 1.979608e+01\\
cvt$_{3}$ & 8192 & 24516 & 1.740063e-02 & 7.226983e-04 & 3.377280 & 2.002835e+01\\
cvt$_{4}$ & 16384 & 49019 & 1.211312e-02 & 5.040025e-04 & 3.345637 & 2.174663e+01\\
cvt$_{5}$ & 32768 & 98064 & 8.695573e-03 & 3.402956e-04 & 3.449458 & 2.140898e+01\\
cvt$_{6}$ & 65536 & 196067 & 6.110583e-03 & 2.334294e-04 & 3.675635 & 2.278220e+01\\
\bottomrule
\end{tabular}
\end{table}

\

\begin{table}
\centering
\caption{Random Voronoi cells used in experiment i).}
\label{tab:mesh-voro}
\begin{tabular}{
c
S[table-format=5.0]
S[table-format=6.0]
S[table-format=1.{\roundPrecision}e-1]
S[table-format=1.{\roundPrecision}e-1]
S[table-format=1.{\roundPrecision}e+1]
S[table-format=1.{\roundPrecision}e+1]
}
\toprule
Mesh & {$N_p$} & {$N_e$} & {$h$} & {$h_\textup{min}$} & {$\gamma_0$} & {$\gamma_1$}\\
\midrule
voro$_{1}$ & 2500 & 7505 & 6.384666e-02 & 6.192528e-06 & 1.422215e+01 & 5.852731e+03\\
voro$_{2}$ & 5000 & 15007 & 4.344562e-02 & 5.845179e-07 & 1.456663e+01 & 3.404522e+04\\
voro$_{3}$ & 10000 & 30006 & 3.470002e-02 & 1.732139e-07 & 2.525383e+01 & 9.493668e+04\\
voro$_{4}$ & 20000 & 60010 & 2.405393e-02 & 2.138871e-07 & 2.087832e+01 & 7.246942e+04\\
voro$_{5}$ & 40000 & 120006 & 1.726980e-02 & 8.256465e-08 & 2.675024e+01 & 7.178744e+04\\
voro$_{6}$ & 80000 & 240027 & 1.140086e-02 & 5.998477e-09 & 2.881822e+01 & 1.100228e+06\\
\bottomrule
\end{tabular}
\end{table}

The goal of this section is to discuss in greater detail the numerical implementation of our method and to provide evidence of the theoretical estimates proven in Section~\ref{sect:theory}.

As a basis for $\mathbb P_k(K)$, for each $K\in\mathcal T_h$, we use the scaled monomials of degree less then or equal to $k$
\begin{equation}\label{eq:monomials}
m_{\boldsymbol\alpha}(x,y) = \left(\frac{x-x_K}{h_K}\right)^{\alpha_1}\left(\frac{y-y_K}{h_K}\right)^{\alpha_2}
\quad\forall\,\boldsymbol{\alpha} = (\alpha_1,\alpha_2)\in\mathbb N^2,\quad\alpha_1+\alpha_2 \leq k,
\end{equation}
where $(x_K, y_K)$ are the coordinates of the barycenter of $K$. Moreover, as a basis for $\mathbb P_k(e)$, for each $e\in\mathcal E_h$, we use Legendre polynomials of degree $\leq k$.

In order to construct the stabilization form described in Section~\ref{sect:stab}, we let $\hat T$ be the unit triangle $\Set{(x,x)\in\mathbb R^2 | 0\leq x,y\leq 1, x+y \leq 1}$ and $\widehat V_\delta$ be the conforming finite element space of polynomial order $1$ constructed on a mesh $\hat T_h$ of $\hat T$ of mesh size $\delta = (k')^{-2}$, whose elements vanish on $\partial\hat T\setminus\hat e$
\[
\widehat V_\delta = \Set{v \in C^0(\overline{\hat T}) | v|_E\in \mathbb P_1(E)\quad\forall E\in\hat T_h, v|_{\partial\hat T\setminus\hat e} = 0}.
\]
Note that the results shown in the following suggest that the choice $\delta = (k')^{-2}$
may be too conservative. With this definition of $\hat V_\delta$, a basis for the space $\widehat W \subseteq H^1(\hat T)$ is built once for all during the pre-processing phase of Problem~\ref{PbGlob} by solving~\ref{eq:defphihat}, for $\lambda = \lambda_i, i = 0, \dots, k$, being $\lambda_i$ a basis of $\mathbb P_k(\hat e)$. Then, a basis $\{\varphi_i\}_{i=1}^{n_K}$ for the auxiliary space $W^K \subseteq H^1(K)$ is computed as indicated by equations~\eqref{eq:defWi} and \eqref{eq:defWK}. In order to assemble the stabilization term on $K$, the next steps are:
\begin{enumerate}
    \item Assemble the stiffness matrix $S$ associated with $\varphi_i\in W^K$, i.e., \[S_{i,j} = (\nabla\varphi_j,\nabla\varphi_i)_K,\quad\text{ for }i, j = 1, \dots, n_K.\]
    \item For $u, v \in \mathbb P_k(K), \lambda, \mu \in \Lambda^K, \varphi_i \in W^K$, compute $\eta_i, \zeta_i$, and $f_i$ from~\ref{eq:sK-howto-vectors}.\label{enum:solveS}
    \item Solve $S\vec\gamma = \vec\zeta$, and compute $\vec \eta^T\vec \gamma, \vec f^T\vec \gamma$ from~\ref{eq:sK-howto}.
\end{enumerate}
Since functions in $W_i \subset W_K$ have zero Dirichlet traces along the edges interior to $K$, the stiffness matrix $S$ is block diagonal, with blocks of size $(k+1)\times(k+1)$, since dim$(W_i) = k+1$, thereby decoupling the contribution of the $n_K$ triangles $T_i$ to the stabilization term $s_K  $. Thus, at step~\ref{enum:solveS}, one has to solve $n_K$ small systems of dimension $k+1$, rather than a single big system of dimension $n_K\cdot(k+1)\times n_K\cdot(k+1)$. Moreover, for each $T_i$, computing the terms
\begin{align*}
&S|_{T_i}, & &\int_{T_i}\nabla v\cdot\nabla\varphi_i, & &\int_{\partial T_i\cap\partial K}\mu\,\varphi_i, & &\int_{T_i}f\varphi_i
\end{align*}
coming from~\eqref{eq:sK-howto-vectors} does not require much effort: one can store the stiffness matrix, the right hand side, and the nodal values of the basis $\varphi_i$ computed only once during the pre-processing phase, and then apply proper push-back operations between $T_i$ and $\hat T$, which amount to matrix-matrix, or matrix-vector multiplications, efficiently performed in our code using PETSc interfaces to BLAS/LAPACK software~\cite{petsc-user-ref}.

After computing the stabilization term, locally for each $K$, problem~\ref{PbGlob} is solved using \emph{static condensation}: for each $K$, equation \eqref{localpb}  yields a local discrete Dirichlet problem, thereby allowing to express $u|_K, \lambda|_{\partial K}$ as a function of the sole variable $\varphi|_{\partial K}$. At this point, we use \eqref{eq:PbGlob-glue}, which imposes continuity of the fluxes $\lambda$, to glue all the local problems together and obtain a global system of equations where only $\varphi$ appears as unknown. The global system is solved with the direct solver STRUMPACK~\cite{strumpack}. Reconstruction of $u, \lambda$ is done by solving local problems in parallel.

We performed a series of experiments in order to investigate the performance of our method with regards to: i) optimal order of convergence of $\|u - u_h\|_{H^1}$; ii) robustness for increasing polynomial degree $k$; iii) sensitivity to the choice of the mesh size $\delta$; iv) robustness with respect to collapsing minimum edge length.

In all the experiments, we let the domain $\Omega$ be the unit square $[0,1]\times[0,1]$. Problem \ref{PbGlob} is solved with Neumann boundary conditions on $\Gamma_N = \{(x,y)\,|\,0 \leq x \leq 1, y = 1\}$, Dirichlet boundary conditions on $\Gamma_D = \partial\Omega\setminus\Gamma_N$, and load term chosen in such a way that
\[
u = \frac{1}{128\pi^2}\cos(8\pi x)\cos(8\pi y)
\]
is the exact solution. The stabilization parameters are chosen to be $\alpha = t = 1$. For the first three experiments, we consider three types of meshes: meshes made mainly of regular hexagons (see, e.g., Figure~\ref{fig:rhexa}), Central Voronoi Tessellations (see, e.g., Figure~\ref{fig:cvt}), and random Voronoi meshes (see, e.g., Figure~\ref{fig:voro}). Geometrical data for the meshes used in experiment i) are shown in Tables~\ref{tab:mesh-r-hexa}, \ref{tab:mesh-cvt}, \ref{tab:mesh-voro}, respectively. For each mesh, we provide: $N_p$, the number of elements of $\mathcal T_h$; $N_e$, the number of edges of $\mathcal T_h$; $h = \max_{K\in\Omega_h} h_K$; $h_\textup{min} = \min_{K\in\mathcal T_h} h_{\textup{min},K}$, where $h_{\textup{min},K}$ is the minimum distance between any two vertices of $K$; $\gamma_0 = \max_{K\in\mathcal T_h}\frac{h_K}{\rho_K}$, where $\rho_K$ is the radius of the largest ball that is contained inside $K$; $\gamma_1 = \max_{K\in\mathcal T_h}\frac{h_K}{h_{\textup{min},K}}$.

\begin{table}
\centering
\caption{Errors and estimated convergence rates (ecr) for experiment i) on random Voronoi cells, $k = 1,2$.}
\begin{tabular}{
S[table-format=6.0]
S[table-format=1.{\roundPrecision}e-1]
S[table-format=1.2]
S[table-format=1.{\roundPrecision}e-1]
S[table-format=1.2]
|
S[table-format=6.0]
S[table-format=1.{\roundPrecision}e-1]
S[table-format=1.2]
S[table-format=1.{\roundPrecision}e-1]
S[table-format=1.2]
}
\toprule
\multicolumn{5}{c}{{$k = 1$}} & \multicolumn{5}{c}{{$k = 2$}} \\
\midrule
{dofs} & {$e^u_1$} & {ecr} & {$e^u_0$} & {ecr} & {dofs} & {$e^u_1$} & {ecr} & {$e^u_0$} & {ecr}\\
\midrule
7500   &   2.616061e-01   &   {-}   &   8.573415e-02   &   {-}   &   15000   &   3.837797e-02   &   {-}   &   1.080663e-02   &   {-}\\
15000   &   1.866574e-01   &   0.974007   &   4.481604e-02   &   1.871709   &   30000   &   1.930107e-02   &   1.983196   &   4.022034e-03   &   2.851839\\
30000   &   1.300656e-01   &   1.042306   &   2.292847e-02   &   1.933750   &   60000   &   9.664237e-03   &   1.995906   &   1.523124e-03   &   2.801784\\
60000   &   9.189962e-02   &   1.002218   &   1.131450e-02   &   2.037934   &   120000   &   4.866747e-03   &   1.979396   &   6.370479e-04   &   2.515119\\
120000   &   6.488821e-02   &   1.004205   &   5.702308e-03   &   1.977111   &   240000   &   2.457406e-03   &   1.971643   &   2.909055e-04   &   2.261703\\
240000   &   4.594882e-02   &   0.995857   &   2.934847e-03   &   1.916521   &   480000   &   1.220759e-03   &   2.018715   &   1.297949e-04   &   2.328633\\
\bottomrule
\end{tabular}
\label{tab:h_convergence-voro-k12}
\end{table}

\begin{table}
\centering
\caption{Errors and estimated convergence rates (ecr) for experiment i) on random Voronoi cells, $k = 3,4$.}
\begin{tabular}{
S[table-format=6.0]
S[table-format=1.{\roundPrecision}e-1]
S[table-format=1.2]
S[table-format=1.{\roundPrecision}e-1]
S[table-format=1.2]
|
S[table-format=6.0]
S[table-format=1.{\roundPrecision}e-1]
S[table-format=1.2]
S[table-format=1.{\roundPrecision}e-1]
S[table-format=1.2]
}
\toprule
\multicolumn{5}{c}{{$k = 3$}} & \multicolumn{5}{c}{{$k = 4$}} \\
\midrule
{dofs} & {$e^u_1$} & {ecr} & {$e^u_0$} & {ecr} & {dofs} & {$e^u_1$} & {ecr} & {$e^u_0$} & {ecr}\\
\midrule
25000   &   3.994555e-03   &   {-}   &   7.052251e-04   &   {-}   &   37500   &   3.690449e-04   &   {-}   &   5.713895e-05   &   {-}\\
50000   &   1.467950e-03   &   2.888463   &   1.911587e-04   &   3.766625   &   75000   &   8.904417e-05   &   4.102406   &   1.011848e-05   &   4.994964\\
100000   &   5.085028e-04   &   3.058951   &   4.637896e-05   &   4.086457   &   150000   &   2.215651e-05   &   4.013582   &   1.876430e-06   &   4.861860\\
200000   &   1.779239e-04   &   3.029990   &   1.164958e-05   &   3.986386   &   300000   &   5.552195e-06   &   3.993200   &   3.863092e-07   &   4.560325\\
400000   &   6.409269e-05   &   2.946058   &   3.042420e-06   &   3.873973   &   600000   &   1.412575e-06   &   3.949462   &   7.707357e-08   &   4.650896\\
800000   &   2.254196e-05   &   3.015093   &   7.652456e-07   &   3.982449   &   1200000   &   3.488335e-07   &   4.035433   &   1.665190e-08   &   4.421098\\
\bottomrule
\end{tabular}
\label{tab:h_convergence-voro-k34}
\end{table}

\begin{table}
\centering
\caption{Errors and estimated convergence rates (ecr) for experiment i) on random Voronoi cells, $k = 5,6$.}
\begin{tabular}{
S[table-format=7.0]
S[table-format=1.{\roundPrecision}e-1]
S[table-format=1.2]
S[table-format=1.{\roundPrecision}e-2]
S[table-format=1.2]
|
S[table-format=7.0]
S[table-format=1.{\roundPrecision}e-2]
S[table-format=1.2]
S[table-format=1.{\roundPrecision}e-2]
S[table-format=1.2]
}
\toprule
\multicolumn{5}{c}{{$k = 5$}} & \multicolumn{5}{c}{{$k = 6$}} \\
\midrule
{dofs} & {$e^u_1$} & {ecr} & {$e^u_0$} & {ecr} & {dofs} & {$e^u_1$} & {ecr} & {$e^u_0$} & {ecr}\\
\midrule
52500   &   2.704856e-05   &   {-}   &   3.316544e-06   &   {-}   &   70000   &   2.183077e-06   &   {-}   &   2.349596e-07   &   {-}\\
105000   &   4.852633e-06   &   4.957424   &   4.166513e-07   &   5.985536   &   140000   &   2.287170e-07   &   6.509456   &   1.828481e-08   &   7.367391\\
210000   &   8.344122e-07   &   5.079871   &   5.267484e-08   &   5.967309   &   280000   &   2.891343e-08   &   5.967503   &   1.684279e-09   &   6.880884\\
420000   &   1.413387e-07   &   5.123208   &   6.367020e-09   &   6.096848   &   560000   &   3.655275e-09   &   5.967376   &   1.666250e-10   &   6.674909\\
840000   &   2.660439e-08   &   4.818841   &   8.687998e-10   &   5.747045   &   1120000   &   4.580953e-10   &   5.992520   &   1.634106e-11   &   6.700064\\
1680000   &   4.555483e-09   &   5.091976   &   1.022008e-10   &   6.175234   &   2240000   &   5.679056e-11   &   6.023849   &   1.664245e-12   &   6.591124\\
\bottomrule
\end{tabular}
\label{tab:h_convergence-voro-k56}
\end{table}

%%%%%%%%%%%%%%%%%%%%%%%%%%%%%%%%%%%%%%

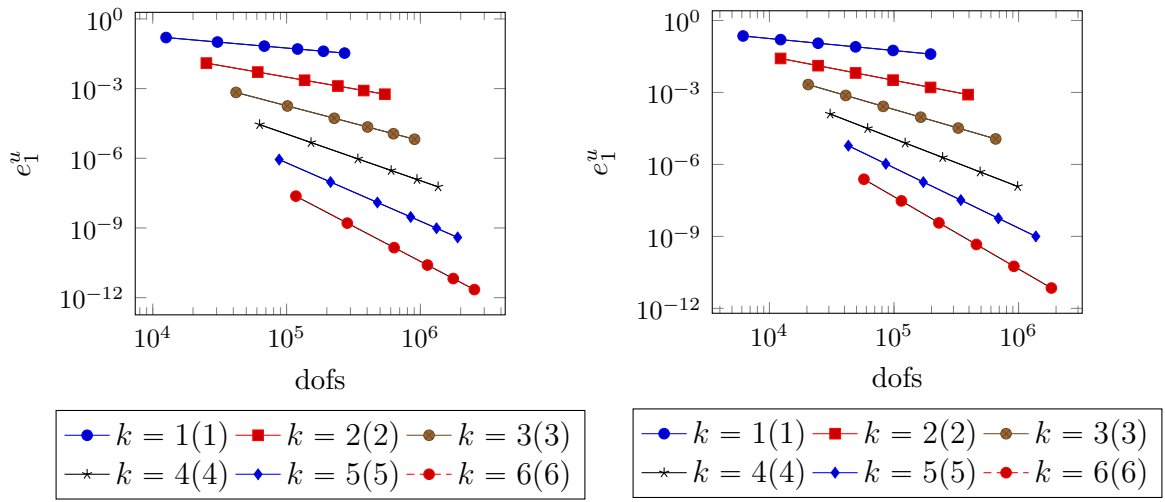
\begin{figure}
\centering
\subfloat{
\begin{tabular}{c}%{rl}
\begin{tikzpicture}[trim axis left]
\begin{loglogaxis}[small, xlabel={dofs}, ylabel={$e^u_1$},
legend columns=3, legend to name=h-rhexa-H1, cycle list name=color]
\addplot[sharp plot,forget plot] table[row sep=\\,y={create col/linear regression={y=err}}]
{
dofs err\\
12579 0.157399\\
30453 0.100879\\
68178 0.0673059\\
120903 0.0504997\\
188628 0.0404097\\
271353 0.0336804\\
};
\expandafter\xdef\csname H1-slope-1\endcsname{\pgfplotstableregressiona}
\addplot coordinates {
(12579,0.157399)
(30453,0.100879)
(68178,0.0673059)
(120903,0.0504997)
(188628,0.0404097)
(271353,0.0336804)
};
\addlegendentry{$k = 1 (\pgfmathparse{(\csname H1-slope-1\endcsname)*(-2)}\pgfmathprintnumber[precision=1,fixed]{\pgfmathresult})$}
\addplot[sharp plot,forget plot] table[row sep=\\,y={create col/linear regression={y=err}}]
{
dofs err\\
25158 0.0124913\\
60906 0.00513293\\
136356 0.00228513\\
241806 0.0012864\\
377256 0.000823668\\
542706 0.000572161\\
};
\expandafter\xdef\csname H1-slope-2\endcsname{\pgfplotstableregressiona}
\addplot coordinates {
(25158,0.0124913)
(60906,0.00513293)
(136356,0.00228513)
(241806,0.0012864)
(377256,0.000823668)
(542706,0.000572161)
};
\addlegendentry{$k = 2 (\pgfmathparse{(\csname H1-slope-2\endcsname)*(-2)}\pgfmathprintnumber[precision=1,fixed]{\pgfmathresult})$}
\addplot[sharp plot,forget plot] table[row sep=\\,y={create col/linear regression={y=err}}]
{
dofs err\\
41930 0.000683456\\
101510 0.000178357\\
227260 5.2746e-05\\
403010 2.22395e-05\\
628760 1.13844e-05\\
904510 6.58773e-06\\
};
\expandafter\xdef\csname H1-slope-3\endcsname{\pgfplotstableregressiona}
\addplot coordinates {
(41930,0.000683456)
(101510,0.000178357)
(227260,5.2746e-05)
(403010,2.22395e-05)
(628760,1.13844e-05)
(904510,6.58773e-06)
};
\addlegendentry{$k = 3 (\pgfmathparse{(\csname H1-slope-3\endcsname)*(-2)}\pgfmathprintnumber[precision=1,fixed]{\pgfmathresult})$}
\addplot[sharp plot,forget plot] table[row sep=\\,y={create col/linear regression={y=err}}]
{
dofs err\\
62895 2.82408e-05\\
152265 4.76106e-06\\
340890 9.42828e-07\\
604515 2.98666e-07\\
943140 1.22416e-07\\
1.35676e+06 5.90616e-08\\
};
\expandafter\xdef\csname H1-slope-4\endcsname{\pgfplotstableregressiona}
\addplot coordinates {
(62895,2.82408e-05)
(152265,4.76106e-06)
(340890,9.42828e-07)
(604515,2.98666e-07)
(943140,1.22416e-07)
(1.35676e+06,5.90616e-08)
};
\addlegendentry{$k = 4 (\pgfmathparse{(\csname H1-slope-4\endcsname)*(-2)}\pgfmathprintnumber[precision=1,fixed]{\pgfmathresult})$}
\addplot[sharp plot,forget plot] table[row sep=\\,y={create col/linear regression={y=err}}]
{
dofs err\\
88053 8.79521e-07\\
213171 9.49057e-08\\
477246 1.25307e-08\\
846321 2.97723e-09\\
1.3204e+06 9.76272e-10\\
1.89947e+06 3.92523e-10\\
};
\expandafter\xdef\csname H1-slope-5\endcsname{\pgfplotstableregressiona}
\addplot coordinates {
(88053,8.79521e-07)
(213171,9.49057e-08)
(477246,1.25307e-08)
(846321,2.97723e-09)
(1.3204e+06,9.76272e-10)
(1.89947e+06,3.92523e-10)
};
\addlegendentry{$k = 5 (\pgfmathparse{(\csname H1-slope-5\endcsname)*(-2)}\pgfmathprintnumber[precision=1,fixed]{\pgfmathresult})$}
\addplot[sharp plot,forget plot] table[row sep=\\,y={create col/linear regression={y=err}}]
{
dofs err\\
117404 2.34841e-08\\
284228 1.61824e-09\\
636328 1.423e-10\\
1.12843e+06 2.53463e-11\\
1.76053e+06 6.64476e-12\\
2.53263e+06 2.22605e-12\\
};
\expandafter\xdef\csname H1-slope-6\endcsname{\pgfplotstableregressiona}
\addplot coordinates {
(117404,2.34841e-08)
(284228,1.61824e-09)
(636328,1.423e-10)
(1.12843e+06,2.53463e-11)
(1.76053e+06,6.64476e-12)
(2.53263e+06,2.22605e-12)
};
\addlegendentry{$k = 6 (\pgfmathparse{(\csname H1-slope-6\endcsname)*(-2)}\pgfmathprintnumber[precision=1,fixed]{\pgfmathresult})$}
\end{loglogaxis}
\end{tikzpicture}
\\
\ref{h-rhexa-H1}
\end{tabular}
}
\subfloat{
\begin{tabular}{c}%{rl}
\begin{tikzpicture}[trim axis left]
\begin{loglogaxis}[small, xlabel={dofs}, ylabel={$e^u_1$},
legend columns=3, legend to name=h-cvt-H1, cycle list name=color]
\addplot[sharp plot,forget plot] table[row sep=\\,y={create col/linear regression={y=err}}]
{
dofs err\\
6144 0.225388\\
12288 0.159138\\
24576 0.112408\\
49152 0.0794934\\
98304 0.0561935\\
196608 0.0397204\\
};
\expandafter\xdef\csname H1-slope-1\endcsname{\pgfplotstableregressiona}
\addplot coordinates {
(6144,0.225388)
(12288,0.159138)
(24576,0.112408)
(49152,0.0794934)
(98304,0.0561935)
(196608,0.0397204)
};
\addlegendentry{$k = 1 (\pgfmathparse{(\csname H1-slope-1\endcsname)*(-2)}\pgfmathprintnumber[precision=1,fixed]{\pgfmathresult})$}
\addplot[sharp plot,forget plot] table[row sep=\\,y={create col/linear regression={y=err}}]
{
dofs err\\
12288 0.0260795\\
24576 0.0130379\\
49152 0.00648373\\
98304 0.00324895\\
196608 0.00161924\\
393216 0.000808493\\
};
\expandafter\xdef\csname H1-slope-2\endcsname{\pgfplotstableregressiona}
\addplot coordinates {
(12288,0.0260795)
(24576,0.0130379)
(49152,0.00648373)
(98304,0.00324895)
(196608,0.00161924)
(393216,0.000808493)
};
\addlegendentry{$k = 2 (\pgfmathparse{(\csname H1-slope-2\endcsname)*(-2)}\pgfmathprintnumber[precision=1,fixed]{\pgfmathresult})$}
\addplot[sharp plot,forget plot] table[row sep=\\,y={create col/linear regression={y=err}}]
{
dofs err\\
20480 0.00212599\\
40960 0.000747797\\
81920 0.000262229\\
163840 9.28492e-05\\
327680 3.26832e-05\\
655360 1.15228e-05\\
};
\expandafter\xdef\csname H1-slope-3\endcsname{\pgfplotstableregressiona}
\addplot coordinates {
(20480,0.00212599)
(40960,0.000747797)
(81920,0.000262229)
(163840,9.28492e-05)
(327680,3.26832e-05)
(655360,1.15228e-05)
};
\addlegendentry{$k = 3 (\pgfmathparse{(\csname H1-slope-3\endcsname)*(-2)}\pgfmathprintnumber[precision=1,fixed]{\pgfmathresult})$}
\addplot[sharp plot,forget plot] table[row sep=\\,y={create col/linear regression={y=err}}]
{
dofs err\\
30720 0.000126284\\
61440 3.15305e-05\\
122880 7.78097e-06\\
245760 1.95194e-06\\
491520 4.84572e-07\\
983040 1.20846e-07\\
};
\expandafter\xdef\csname H1-slope-4\endcsname{\pgfplotstableregressiona}
\addplot coordinates {
(30720,0.000126284)
(61440,3.15305e-05)
(122880,7.78097e-06)
(245760,1.95194e-06)
(491520,4.84572e-07)
(983040,1.20846e-07)
};
\addlegendentry{$k = 4 (\pgfmathparse{(\csname H1-slope-4\endcsname)*(-2)}\pgfmathprintnumber[precision=1,fixed]{\pgfmathresult})$}
\addplot[sharp plot,forget plot] table[row sep=\\,y={create col/linear regression={y=err}}]
{
dofs err\\
43008 6.02437e-06\\
86016 1.05919e-06\\
172032 1.83087e-07\\
344064 3.2513e-08\\
688128 5.68213e-09\\
1.37626e+06 9.98403e-10\\
};
\expandafter\xdef\csname H1-slope-5\endcsname{\pgfplotstableregressiona}
\addplot coordinates {
(43008,6.02437e-06)
(86016,1.05919e-06)
(172032,1.83087e-07)
(344064,3.2513e-08)
(688128,5.68213e-09)
(1.37626e+06,9.98403e-10)
};
\addlegendentry{$k = 5 (\pgfmathparse{(\csname H1-slope-5\endcsname)*(-2)}\pgfmathprintnumber[precision=1,fixed]{\pgfmathresult})$}
\addplot[sharp plot,forget plot] table[row sep=\\,y={create col/linear regression={y=err}}]
{
dofs err\\
57344 2.40875e-07\\
114688 3.01395e-08\\
229376 3.64961e-09\\
458752 4.59313e-10\\
917504 5.64763e-11\\
1.83501e+06 7.03596e-12\\
};
\expandafter\xdef\csname H1-slope-6\endcsname{\pgfplotstableregressiona}
\addplot coordinates {
(57344,2.40875e-07)
(114688,3.01395e-08)
(229376,3.64961e-09)
(458752,4.59313e-10)
(917504,5.64763e-11)
(1.83501e+06,7.03596e-12)
};
\addlegendentry{$k = 6 (\pgfmathparse{(\csname H1-slope-6\endcsname)*(-2)}\pgfmathprintnumber[precision=1,fixed]{\pgfmathresult})$}
\end{loglogaxis}
\end{tikzpicture}
\\
\ref{h-cvt-H1}
\end{tabular}
}
\caption{Experiment i): relative errors $e^u_1$ and convergence rates $(\cdot)$ on hexagonal meshes (left) and CVT meshes (right).}
\label{fig:h-rhexa-cvt}
\end{figure}

\begin{table}
\centering
\caption{History of convergence for increasing polynomial degrees.}
\begin{tabular}{
c
|
S[table-format=1.{\roundPrecision}e-2]
S[table-format=1.2]
|
S[table-format=1.{\roundPrecision}e-2]
S[table-format=1.2]
|
S[table-format=1.{\roundPrecision}e-2]
S[table-format=1.2]
}
\toprule
\multirow{2}*{$k$} & \multicolumn{2}{c}{{r-hexa}} & \multicolumn{2}{c}{{cvt}} & \multicolumn{2}{c}{{voro}}\\
\cmidrule(lr){2-7}
& {$e^u_1$} & {ecr} & {$e^u_1$} & {ecr} & {$e^u_1$} & {ecr}\\
\midrule
$1$   &   9.026205e-01   &   {-}   &   9.757498e-01   &   {-}   &   9.834298e-01   &   {-}\\
$2$   &   6.283715e-01   &   {-}   &   6.314695e-01   &   {-}   &   8.372923e-01   &   {-}\\
$3$   &   2.406498e-01   &   0.377344   &   3.252858e-01   &   0.656003   &   5.406554e-01   &   0.367801\\
$4$   &   8.820264e-02   &   0.956245   &   1.008617e-01   &   0.566500   &   3.345975e-01   &   0.911509\\
$5$   &   2.017019e-02   &   0.680280   &   3.612131e-02   &   1.140317   &   1.654520e-01   &   0.681371\\
$6$   &   5.247257e-03   &   1.095753   &   6.857149e-03   &   0.618002   &   7.854378e-02   &   0.945267\\
$7$   &   8.080464e-04   &   0.719730   &   1.880330e-03   &   1.284228   &   2.993953e-02   &   0.772466\\
$8$   &   1.573773e-04   &   1.143565   &   2.555216e-04   &   0.648253   &   1.095438e-02   &   0.959258\\
$9$   &   1.868786e-05   &   0.767784   &   5.978951e-05   &   1.374133   &   3.262526e-03   &   0.830093\\
$10$   &   2.988154e-06   &   1.162312   &   6.048305e-06   &   0.633978   &   9.631321e-04   &   0.992764\\
$11$   &   2.916116e-07   &   0.787807   &   1.238132e-06   &   1.444389   &   2.338275e-04   &   0.861868\\
$12$   &   3.849414e-08   &   1.149177   &   9.987175e-08   &   0.630066   &   5.753514e-05   &   1.009570\\
$13$   &   3.100948e-09   &   0.803922   &   1.860324e-08   &   1.498004   &   1.178006e-05   &   0.884110\\
$14$   &   3.644692e-10   &   1.176447   &   1.229514e-09   &   0.618596   &   2.494517e-06   &   1.021693\\
$15$   &   2.708842e-10   &   7.214883   &   2.828873e-10   &   1.848955   &   4.996221e-06   &   -2.234873\\
$16$   &   1.229638e-09   &   -0.196162   &   6.401413e-10   &   -1.799232   &   1.565970e-06   &   -0.598691\\
\bottomrule
\end{tabular}
\label{tab:k-robustness}
\end{table}

\begin{table}
\centering
\caption{Experiment iii): relative errors $e^u_1$ on the hexagonal mesh shown in Figure~\ref{fig:meshes-k-robustness} for different choices of the mesh size $\delta$.}
\label{tab:delta_robustness}
\begin{tabular}{
S[table-format=2.0]
S[table-format=1.{\roundPrecision}e-2]
S[table-format=1.{\roundPrecision}e-2]
S[table-format=1.{\roundPrecision}e+1]
S[table-format=1.{\roundPrecision}e-2]
}
\toprule
{$k$} & {$\delta = k^{-2}$} & {$\delta = k^{-1}$} & {$\delta = 1/4$} & {$\delta = 1/8$}\\
\midrule
2 &  6.283715e-01 & 1.264991e+02 & 6.283715e-01 & 6.275545e-01\\
3 & 2.406498e-01 & 7.769248e-01 & 1.684592e+01 & 2.406815e-01\\
4 & 8.820264e-02 & 1.300776e+01 & 1.300776e+01 & 8.797154e-02\\
5 & 2.017019e-02 & 1.999461e-02 & 4.158245e+01 & 2.008178e-02\\
6 & 5.247257e-03 & 5.123302e-03 & 5.146298e+01 & 5.130301e-03\\
7 & 8.080464e-04 & 8.771726e-04 & 5.717572e+02 & 8.346685e-04\\
8 & 1.573773e-04 & 1.589790e-04 & 7.967835e+01 & 1.589790e-04\\
9 & 1.868786e-05 & 2.518937e-05 & 3.525097e+01 & 2.517245e-05\\
10 & 2.988154e-06 & 2.957016e-06 & 6.147804e+02 & 3.185149e-06\\
11 & 2.916116e-07 & 2.997204e-07 & 2.180339e+02 & 5.773862e-07\\
12 & 3.849414e-08 & 3.801017e-08 & 7.752028e+03 & 3.986904e-08\\
13 & 3.100948e-09 & 3.532268e-09 & 4.299032e+02 & 9.339793e-09\\
14 & 3.644692e-10 & 3.605954e-10 & 2.606027e+03 & 1.235350e-08\\
15 & 2.708842e-10 & 2.142945e-10 & 2.986630e+04 & 3.247526e-08\\
16 & 1.229638e-09 & 9.598226e-10 & 3.775751e+05 & 4.032767e-07\\
\bottomrule
\end{tabular}
\end{table}

\begin{enumerate}
    \item[i)] \emph{Optimal order of convergence in $h$}: Tables \ref{tab:h_convergence-voro-k12}--\ref{tab:h_convergence-voro-k56} show the relative errors $e^u_1 = ||u-u_h ||_{H^1} / || u ||_{H^1}, e^u_0 = ||u-u_h ||_{L^2} / || u ||_{L^2}$ and the estimated convergence rates (ecr) for several values of the polynomial degree $k$ on the random Voronoi cells versus the total number of the degrees of freedom dofs $= \text{dim}(V_h)$ (notice that dofs behaves like $\mathcal O(h^{-2})$). Analogous results are plotted in Figure~\ref{fig:h-rhexa-cvt} for both hexagonal and CVT meshes. We note that the results confirm the theoretical estimate, with the correct order of convergence for the $H^1$ norm of the error, i.e. $\mathcal O(h^k)$, as $h$ tends to zero. The difference in convergence rates between odd and even values of $k$ for the $L^2$ norm of the error is consistent with results obtained for non symmetric interior penalty approximations of linear elliptic problems~\cite{houston2002,BABUSKA1999103}.
    
    \item[ii)] \emph{Validity as a $k$-method}: we test the validity of our method as a $k$-method, by fixing the mesh (one of those depicted in Figure~\ref{fig:meshes-k-robustness}) and increasing $k$ from $1$ to $16$. We compute the relative errors $e^u_1$ as functions of $k$ and check whether the rates
    \[
    \frac{\log(e_k/e_{k+1})}{\log(e_{k+1}/e_{k+2})} \approx 1,
    \]
    as would be expected. Table~\ref{tab:k-robustness} shows that this is indeed the case. The loss of accuracy at high order, i.e. $k = 15, 16$, is most probably a consequence of the ill-conditioning due to the choice of the monomial basis~\eqref{eq:monomials}.
    
    \item[iii)] \emph{Sensitivity with respect to the mesh size $\delta$}: Table~\ref{tab:delta_robustness} shows that taking $\delta = k^{-2}$ is a conservative choice ensuring that the error decreases with increasing $k$. However, more permissive choices, e.g. $\delta = k^{-1}$, might be enough to compute the stabilization, provided that $\delta$ is small enough when $k$ is also small, say $k = 2, 3$. Letting $\delta$ being a constant, even if small, has a detrimental effect for increasing $k$, see columns corresponding to $\delta = 1/4$ and $\delta = 1/8$ in Table~\ref{tab:delta_robustness}.

\begin{figure}
\centering
\subfloat[]{\includegraphics[width=0.24\textwidth]{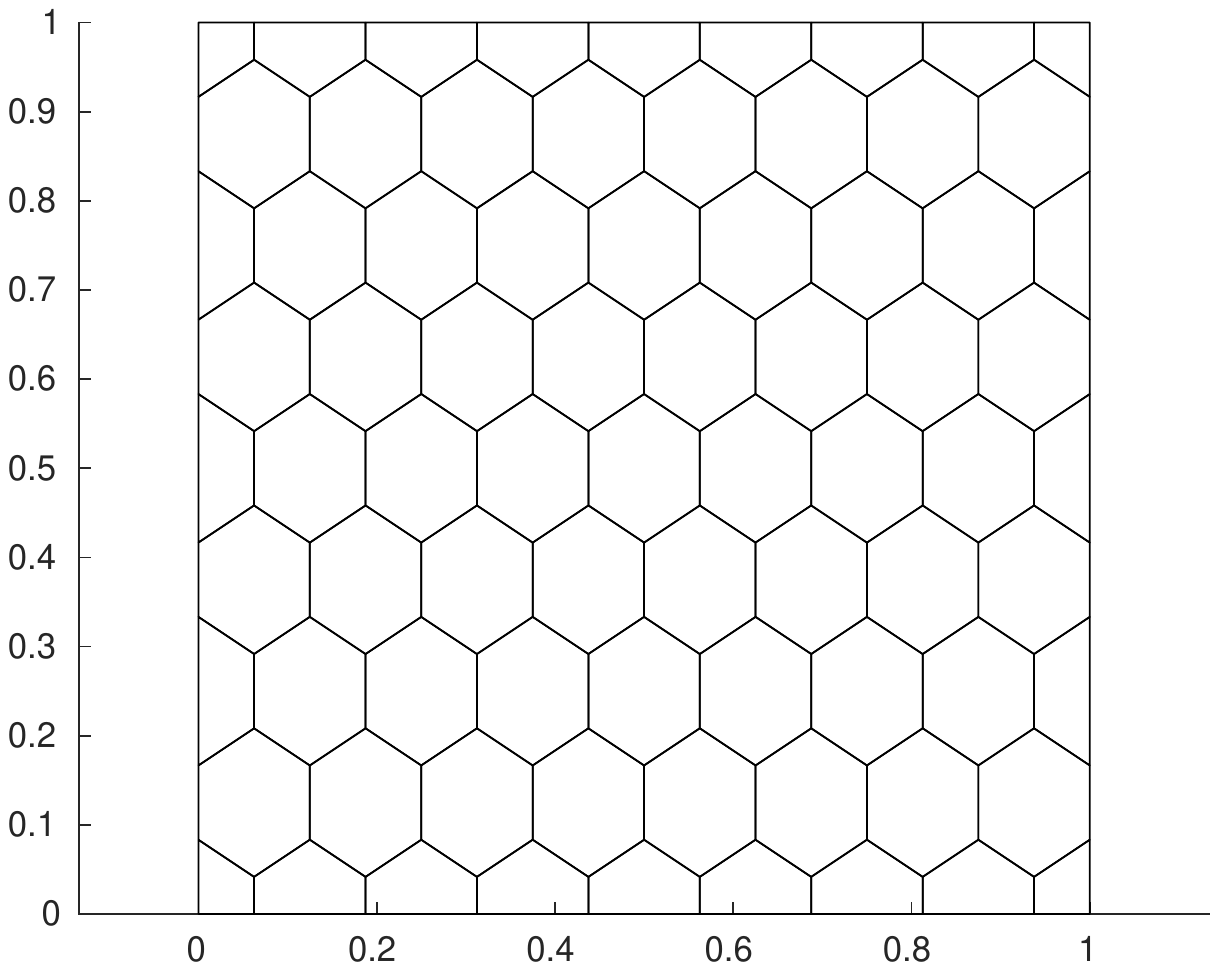}\label{fig:shexa00}}
\subfloat[]{\includegraphics[width=0.24\textwidth]{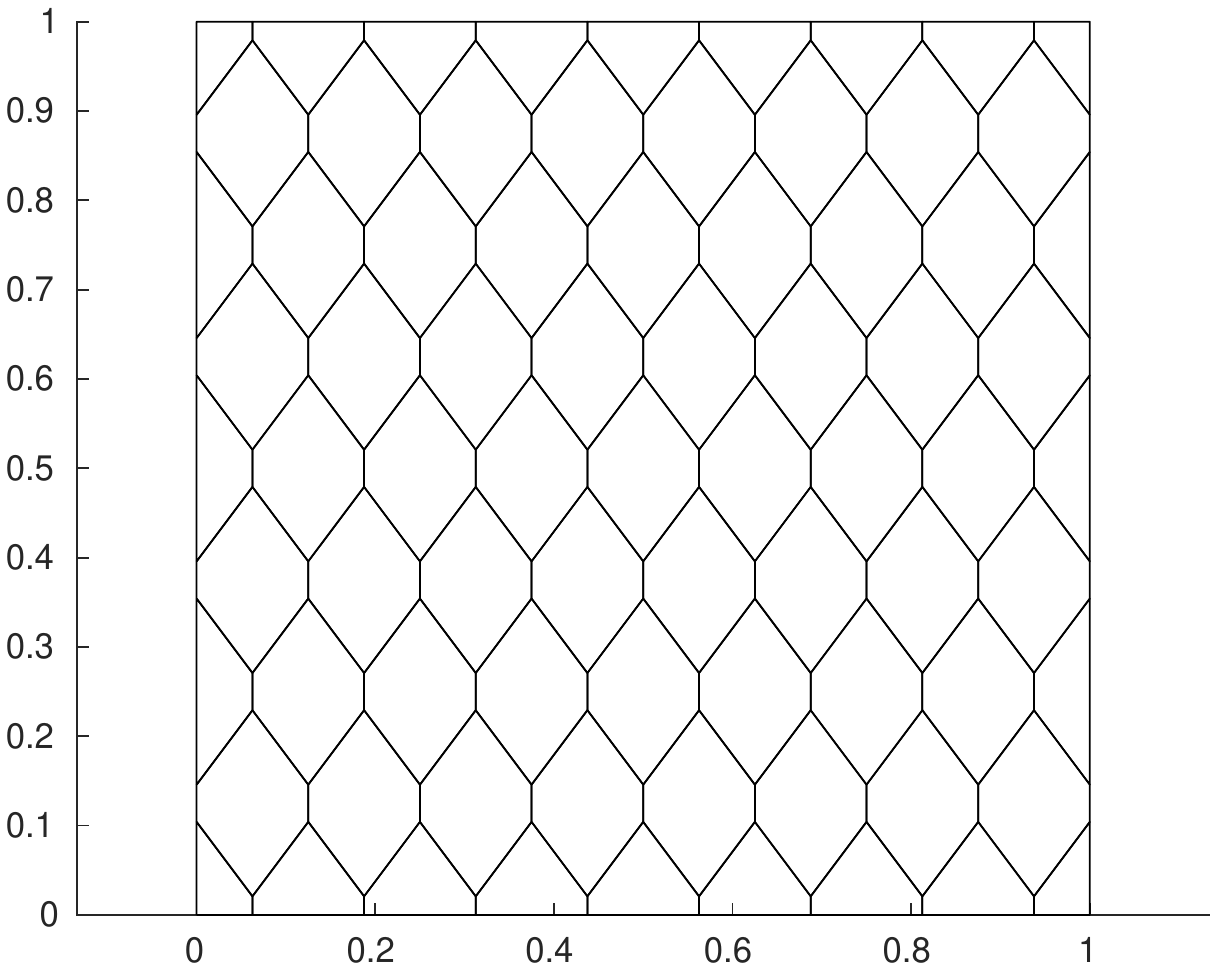}\label{fig:shexa01}}
\subfloat[]{\includegraphics[width=0.24\textwidth]{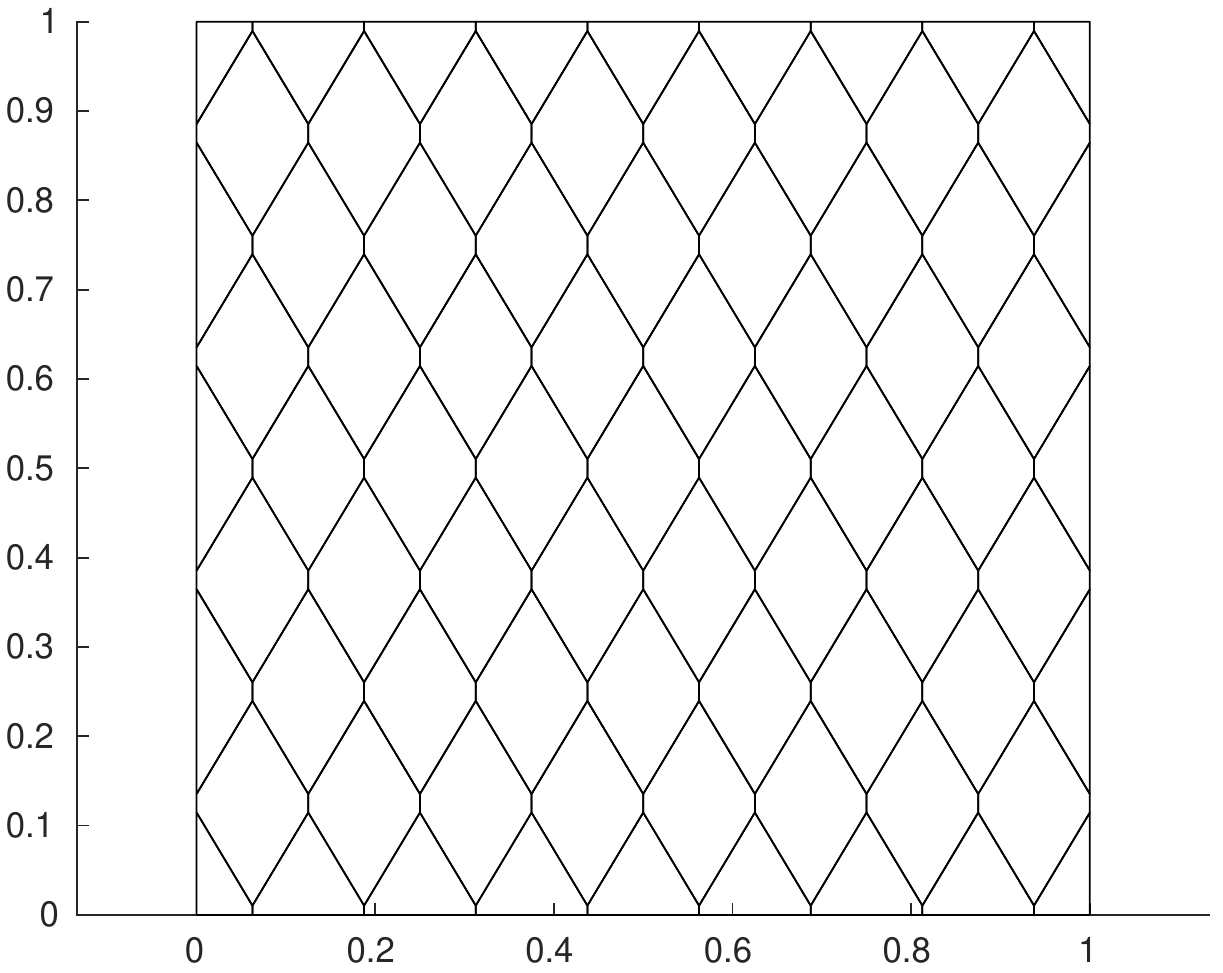}\label{fig:shexa02}}
\subfloat[]{\includegraphics[width=0.24\textwidth]{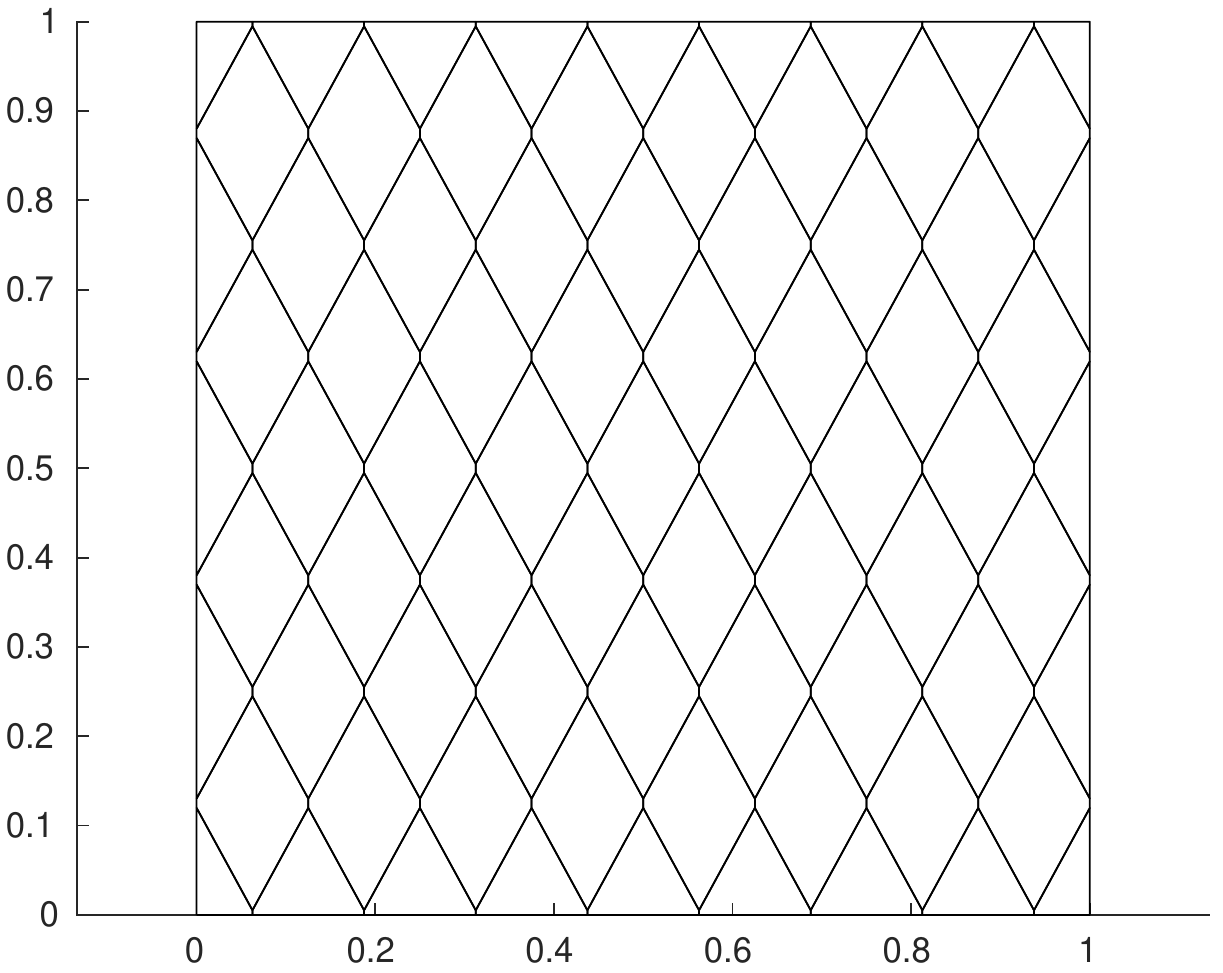}\label{fig:shexa03}}
\caption{Some meshes used in experiment iv). Figure~\ref{fig:shexa00} shows a reference mesh; then, from left to right, we show the meshes obtained by shrinking the vertical edges of a factor $2^{-1}, 2^{-2}, 2^{-3}$, respectively.}
\label{fig:meshes-delta-robustness}
\end{figure}

\begin{table}
\centering
\caption{Experiment iv): History of convergence for shrinking minimum edge length and at different polynomial degrees.}
\label{tab:hmin_robustness}
\begin{tabular}{
c
S[table-format=1.{\roundPrecision}e-2]
S[table-format=1.2]
S[table-format=1.{\roundPrecision}e-2]
S[table-format=1.2]
S[table-format=1.{\roundPrecision}e-2]
S[table-format=1.2]
S[table-format=1.{\roundPrecision}e-2]
S[table-format=1.2]
S[table-format=1.{\roundPrecision}e-1]
S[table-format=1.2]
}
\toprule
\multirow{2}*{$k$} & \multicolumn{2}{c}{{$s = 1$}} & \multicolumn{2}{c}{{$s = 2^{-4}$}} & \multicolumn{2}{c}{{$s = 2^{-8}$}} & \multicolumn{2}{c}{{$s = 2^{-16}$}} & \multicolumn{2}{c}{{$s = 2^{-32}$}}\\
\cmidrule{2-11}
& {$e^u_1$} & {ecr} & {$e^u_1$} & {ecr} & {$e^u_1$} & {ecr} & {$e^u_1$} & {ecr} & {$e^u_1$} & {ecr}\\
\midrule
\multirow{3}*{$6$} & 5.247257e-03   &   {-}      & 1.143072e-02   &   {-}      & 1.270424e-02   &   {-}      & 1.279321e-02   &   {-}      & 1.279356e-02   & {-}\\
                   & 9.098856e-05   &   6.264307 & 2.098650e-04   &   6.176041 & 2.350249e-04   &   6.164306 & 2.367921e-04   &   6.163516 & 2.367992e-04   &   6.163511\\
                   & 1.486324e-06   &   6.141583 & 3.455136e-06   &   6.129902 & 3.874832e-06   &   6.127793 & 3.904512e-06   &   6.127584 & 3.906054e-06   &   6.127040\\
\midrule
\multirow{3}*{$7$} & 8.080464e-04   &   {-}      & 2.121082e-03   &   {-}      & 2.402678e-03   &   {-}      & 2.422578e-03   &   {-}      & 2.422657e-03   & {-}\\
                   & 8.193514e-06   &   7.093244 & 2.229711e-05   &   7.037545 & 2.548120e-05   &   7.023910 & 2.570830e-05   &   7.022945 & 2.570911e-05   &   7.022946\\
                   & 6.851922e-08   &   7.141025 & 1.859361e-07   &   7.145243 & 2.125250e-07   &   7.144985 & 2.144178e-07   &   7.144994 & 2.779145e-07   &   6.757855\\
\midrule
\multirow{3}*{$8$} & 1.573773e-04   &   {-}      & 4.812872e-04   &   {-}      & 5.589005e-04   &   {-}      & 5.644955e-04   &   {-}      & 5.645176e-04   & {-}\\
                   & 6.770337e-07   &   8.417881 & 2.126531e-06   &   8.376618 & 2.480760e-06   &   8.369564 & 2.506356e-06   &   8.369095 & 2.506696e-06   &   8.368946\\
                   & 2.729376e-09   &   8.230191 & 8.799755e-09   &   8.191195 & 1.027501e-08   &   8.189825 & 1.038170e-08   &   8.189728 & 7.783211e-08   &   \textbf{5.18}\\
\midrule
\multirow{3}*{$9$} & 1.868786e-05   &   {-}      & 6.925956e-05   &   {-}      & 8.175467e-05   &   {-}      & 8.265652e-05   &   {-}      & 8.266003e-05   & {-}\\
                   & 4.711766e-08   &   9.243344 & 1.775431e-07   &   9.217731 & 2.110553e-07   &   9.206846 & 2.135054e-07   &   9.205964 & 2.151916e-07   &   9.193876\\
                   & 9.809867e-11   &   9.216535 & 3.691790e-10   &   9.218413 & 4.408099e-10   &   9.211808 & 4.460954e-10   &   9.211244 & 7.136329e-08   &   \textbf{1.65}\\
\bottomrule
\end{tabular}
\end{table}

    \item[iv)] \emph{Robustness with respect to collapsing minimum edge length}: for this experiment, we consider a mesh and two non-nested refinements as reference meshes, and then progressively shrink the length of their vertical edges by a factor of $s = 1 (\text{original mesh}),$ $2^{-1}, 2^{-2}, \dots, 2^{-32}$.
      Convergence is severely and abruptly affected only starting with $k = 8$, on the finest mesh, for the smallest shrinking factor $s = 2^{-32}$ ($h_\textup{min} \approx$ \num{2.43e-12}) (see Table~\ref{tab:hmin_robustness}). 
{Although Assumption~\ref{ass:meshes} (ii) is not satisfied, the method seems quite robust with respect to the minimal edge length, at least for low degrees~$k$, in the approximation of $u$. On the other hand, for $k\ge 8$, the loss of robustness could also be caused by round-off errors.}
% This behavior is likely to be due to round-off errors, rather than to stability issues, thereby suggesting that our numerical scheme is also robust with respect to the minimal edge length.
\end{enumerate}

%\section{}

\bibliographystyle{plain}
\bibliography{biblio2}

% -------------------------------------------------------------
\end{document}